%% file: PflWilk.ECDNDR.tex
\documentclass[11pt]{amsart}
\usepackage{amssymb, amsmath, eucal, rotating}
\usepackage{color}
\usepackage{mathabx}
\usepackage{enumerate}
\usepackage{hyperref}

\definecolor{darkgreen}{rgb}{0, 0.45, 0}
\definecolor{lightgreen}{rgb}{0.5, 1, 0.5}

\input xy
 \xyoption{all}

\numberwithin{equation}{section}

\theoremstyle{plain}

\newtheorem{proposition}{Proposition}[section]
\newtheorem{theorem}[proposition]{Theorem}		
\newtheorem{corollary}[proposition]{Corollary}
\newtheorem{lemma}[proposition]{Lemma}

\theoremstyle{definition}

\newtheorem{definition}[proposition]{Definition}
\newtheorem{remark}[proposition]{Remark}
\newtheorem{example}[proposition]{Example}

\newcommand{\R}{\mathbb R}

\newcommand{\PBbb}{\mathbb P}
\newcommand{\NBbb}{\mathbb N}

\newcommand{\calR}{\mathcal{R}}
\newcommand{\calS}{\mathcal{S}}
\newcommand{\calZ}{\mathcal{Z}}
\newcommand{\calT}{\mathcal{T}}
\newcommand{\calV}{\mathcal{V}}
\newcommand{\calH}{\mathcal{H}}
\newcommand{\sfT}{\mathsf{T}}

\newcommand{\depth}{\operatorname{dp}}
\newcommand{\height}{\operatorname{ht}}

\newcommand{\Hom}{\mathop{\rm Hom}\nolimits}

\newcommand{\tr}{\mathop{\rm Tr}\nolimits}

\newcommand{\SU}{\mathsf{SU}}

\DeclareMathOperator{\id}{id}

\DeclareMathOperator{\supp}{supp}

\begin{document}

\title{Equivariant control data and neighborhood deformation retractions}

\author{Markus J. Pflaum}
\address{Department of Mathematics, University of Colorado, Boulder, CO 80309-0395}
\email{markus.pflaum@colorado.edu}

\author{Graeme Wilkin}
\address{Department of Mathematics, National University of Singapore, Singapore 119076}
\email{graeme@nus.edu.sg}

\date{\today}

\begin{abstract} 
In this article we study Whitney (B) regular stratified spaces with the action of a compact Lie group $G$ which preserves the strata. We prove an equivariant submersion theorem and use it to show that such a $G$-stratified space carries a system of $G$-equivariant control data. As an application, we show that if $A \subset X$ is a closed $G$-stratified subspace which is a union of strata of $X$, then the inclusion $i : A \hookrightarrow X$ is a $G$-equivariant cofibration. In particular, this theorem applies whenever $X$ is a $G$-invariant analytic subspace of an analytic $G$-manifold $M$ and $A \hookrightarrow X$ is a closed $G$-invariant analytic subspace of $X$.
\end{abstract}

\subjclass{}

\dedicatory{In memory of John Mather}

\maketitle

%\tableofcontents

\thispagestyle{empty}

\baselineskip=16pt
%\thispagestyle{plain}
%\setcounter{page}{1}

%\setcounter{section}{-1}

%%%%%%%%%%%%%%%%%%%%%%%%%%%%%%%%%%%%%%%%%%%%%%%%

\input Introduction.tex

\input Equivariant_control_data.tex

\input Equivariant_Wall_result.tex

\input Equivariant_deformation_retract.tex

\appendix
\input Appendix_stratifications.tex

\nocite{FieTGM}

\bibliography{lmlib}
%\addcontentsline{toc}{section}{References}
\bibliographystyle{acm}

\end{document}

%% file: Introduction.tex
\section{Introduction}

Mather's concept of control data \cite{MatNTS} has crystallized as an indispensible tool for the proof of Thom's first and second isotopy lemmata
and more generally for the proof of the topological stability theorem
which was originally conjectured by Thom \cite{ThoEMS} and finally proved by Mather \cite{MatNTS} % GibWirPlesLooTSSM}. %Gibson, Wirthmueller, du Plessis, Loojenga, Toplogical Stability of Smooth Mappings
Moreover, control data are a powerful tool in stratified Morse theory 
\cite{GorPheSMT}, to prove triangulability of stratified spaces fulfilling Whitney's condition (B) \cite{GorTSS}, 
and to verify de Rham theorems in intersection homology theory \cite{BraHecSarTRVS}. A further topological application 
of the concept of control data is that it allows for a transparent proof that every submersed stratified subspace $A$ of a (B) 
regular stratified space $X$ is a neighborhood deformation retract (NDR) or equivalently that 
$i: A \hookrightarrow X$ is a cofibration.  The assumption that $A$ is a closed submersed stratified subspace of $X$ 
hereby means that $A$ is a union of connected components of strata of $X$; see Appendix \ref{sec:stratified-spaces}.

In this article we extend the existence of control data and the latter result 
to the $G$-equivariant case, where $G$ is a compact Lie group. 
More precisely, we show in Theorem \ref{thm:control-data-existence} that if $M$ is a smooth $G$-manifold and $X \subset M$
a (B) regular stratified subspace such that $G$ leaves the strata invariant, then there exists a system of 
$G$-equivariant control data on $X$. We use this observation in Section \ref{sec:equivariant-def-retract} to prove that for every $G$-invariant closed 
submersed stratified subspace $A \subset X$ the inclusion $i: A \hookrightarrow X$ is a $G$-cofibration. More precisely, 
we prove the following which is the main result of our paper.

\begin{theorem}
Let $X$ be a $G$-invariant (B) regular stratified space in a $G$-manifold $M$ and $A$ a $G$-invariant 
closed submersed stratified subspace of $X$. Then there exists a $G$-invariant open neighborhood $U$ of $A$ in $X$ and a 
stratified $G$-equivariant strong deformation retraction $r: U \times I \to U$ of $U$ onto $A$ such that 
$U_s := r(U,s)$ for each $s \in I= [0,1] $ satisfies 
\begin{enumerate}[{\rm (a)}]
\item $U_s$ is open in $X$ for all $s \in [0,1)$, and

\item $U_s = \bigcup_{s < t \leq 1} U_t$ and $\overline{U_s} = \bigcap_{0 \leq t < s} U_t$ for all $s \in (0,1)$. 
\end{enumerate}
In particular $i : A \hookrightarrow X$ is a $G$-equivariant cofibration. 
\end{theorem}

We then consider the situation where $X$ and $A$ are $G$-invariant analytic subspaces of an analytic 
$G$-manifold $M$ with $A \subset X$ being a closed subspace. Using methods by Wall \cite{WalRS} we show in Theorem \ref{thm:equivariant-Wall} that  
$X$ possesses a $G$-invariant (B) regular stratification such that $A$ is a union of strata. 
Hence our main result applies to  such a $G$-invariant analytic pair $(X,A)$.

\noindent \emph{Acknowledgements.} M.P.~was partially supported by a Simons Collaboration Grant, award nr.~359389 and an NSF Conference Grant, DMS-1543812. 
  M.P.~acknowledges hospitality by the National University of Singapore and the Max-Planck-Institute for Mathematics in Bonn, Germany. 
  G.W.~was partially supported by grant number R-146-000-200-112 from the National University of Singapore. 
  G.W.~would also like to thank the University of Colorado for their hospitality during this project.

%% file: Equivariant_control_data.tex
\section{Control data  compatible with a group action}
\label{sec:equivariant-control-data}

Mather proved in \cite{MatNTS} that every (B) regular stratified subspace of a smooth manifold carries a 
system of control data. In this section we extend his result to the $G$-equivariant case.
To this end we first introduce $G$-equivariant versions of stratifications, tubular neighborhoods, their isomorphisms and diffeotopies.
Afterwards we prove a $G$-equivariant submersion theorem. This will be used to derive uniqueness and existence results 
for equivariant tubular neighborhoods. These tools then entail the main result of this section. 

\subsection{Equivariant versions of stratifications and tubular neighborhoods}
\begin{definition}
  Suppose that a compact Lie group $G$ acts on the total space $X$ of a 
  stratified space $(X,\calS)$. The stratification $\calS$ is called a 
  $G$-\emph{stratification} or $G$-\emph{invariant} and $(X,\calS)$ a 
  $G$-\emph{stratified space} if for all $g \in G$ and
  $x \in X$ the set germs $g\calS_x$ and $\calS_{gx}$ coincide
  and if for each open neighborhood $U$ with an $\calS$-inducing decomposition 
  $\calZ$ the map from a piece $R \in \calZ$ to $gR$ given by the $g$-action 
  is a  diffeomorphism of smooth manifolds. We also say in this situation that 
  the $G$-action on $X$ is compatible with the stratification. 
\end{definition}

\begin{example}
  The orbit type stratification of a $G$-manifold $M$ is $G$-invariant since
  the group action leaves the orbit types of the points of $M$  invariant.  
\end{example}

\begin{proposition}
   Each stratum of a $G$-stratified space $(X ,\calS)$ 
   is preserved by the $G$-action. Moreover, $G$ acts smoothly on 
   the strata of $X$. 
\end{proposition}

\begin{proof}
  Let $x$ be a point of a stratum $S$ of $X$. Choose a decomposition $\calZ$ of an
  open neighborhood $U$ of $x$ inducing the stratification $\calS$ over $U$. 
  Then $gU$ is an open neighborhood of $gx$, and $g\calZ$ a decomposition of $gU$.
  Moreover, if $y \in U$ and $R_y$ is the piece of $\calZ$ through $y$, then $gR_y$ 
  is the piece of $g\calZ$ through $gy$. Hence $g\calZ$ induces the stratification
  $\calS$ over $gU$. But that means that $gx$ has the same depth as $x$ and that 
  the dimension of the piece in which $gx$ lies has the same dimension as the 
  piece of $\calZ$ through $x$. So $gx$ and $x$ lie in the same stratum. 
  So we have proved that $G$ acts on each stratum of $X$. This action is
  smooth since it is smooth locally by definition. The claim is proved.   
\end{proof}

\begin{definition}
By $G$-equivariant system of control data % or short a system of $G$-control data 
on a stratified space $(X,\calS)$ with a compatible $G$-action  we understand a 
family $\calT = (T_S,\pi_S,\varrho_S )_{S\in \calS}$ of triples called \emph{tubes} consisting for each
$S\in \calS$ of an open neighborhood $T_S$ of $S$, a continuous retraction  $\pi_S :T_S \to S$ called
\emph{projection} and a continuous map $\varrho_S :S \to [0,\infty )$ called
\emph{tubular function} such that the following control conditions hold true:
   
\begin{enumerate}[(CC1)]
\item 
  For any $S \in \calS$ the neighborhood $T_S$ is $G$-invariant, 
  the projection $\pi_S$ is $G$-equivariant, and the tubular function $\varrho_S$ is $G$-invariant. 
\item 
  For any $S \in \calS$ the tubular function $\varrho_S$ satisfies $S = \rho_S^{-1}(0)$.
\item  For any $R,S \in \mathcal{S}$ with $R < S$, the map
\begin{equation*}
  (\pi_{R,S}, \rho_{R,S}) : T_{R,S} \rightarrow R \times (0, \infty)
\end{equation*}
is a smooth submersion, where $T_{R,S} := T_R \cap S$, $\pi_{R,S} := \left. \pi_R \right|_{T_{R,S}}$ and $\rho_{R,S} = \left. \rho_R \right|_{T_{R,S}}$. 
\item
 For any $Q,R,S\in \calS$ with $Q < R < S$ 
\begin{equation*}
\pi_{Q,R} \circ \pi_{R,S}|_{T_{Q,R,S}} = \pi_{Q,S}|_{T_{Q,R,S}} \quad \text{and} \quad 
\rho_{Q,R} \circ \pi_{R,S}|_{T_{Q,R,S}} = \rho_{Q,S}|_{T_{Q,R,S}} ,
\end{equation*}
where $T_{Q,R,S} = \pi_{R,S}^{-1}(T_{Q,R})$.
\end{enumerate}
If $G$ is the trivial group one recovers the original definition of a system of control data by Mather \cite{MatNTS}. 
A stratified space $(X,\calS)$ together with some control data $\calT$ will be called a \emph{Thom-Mather stratified space}.
If  $(X,\calS)$ carries a compatible $G$-action and the control data are $G$-equivariant in the sense defined above,
we call $(X,\calS,\calT)$  a \emph{Thom-Mather stratified space} with a \emph{compatible $G$-action}
or briefly a \emph{Thom-Mather $G$-stratified space}. 
\end{definition}

We will now introduce some further language where $G$ always denotes a compact Lie group and $G$ a smooth
$G$-manifold. Let $S \subset M$ be a $G$-invariant smooth manifold. 
By a $G$-\emph{equivariant tubular neighborhood} of $S$ in $M$ we understand a
triple $\sfT = (E, \varepsilon, \varphi)$ where $\pi_E:E \to S$ is a $G$-vector bundle over $S$ carrying a  
$G$-invariant inner product $\langle - , - \rangle : S \to E\otimes E$, $\varepsilon : S \to (0,\infty)$ is a $G$-invariant smooth map, 
and $\varphi$ is a $G$-equivariant diffeomorphism from $B(\varepsilon, E) := \{ v \in E \mid \langle v,v \rangle < \varepsilon (\pi_E(v)) \}$
onto an open neighborhood $T_S$ of $S$   such that 
the composition of $\varphi$ with the zero section of $E$ coincides with the identical embedding of $S$ into $M$. 
Note that by the requirements $T_S$ is a $G$-invariant open neighborhood of $S$ in $M$. 
% It is called the \emph{total space} of the tubular neighborhood $\sfT$ or sometimes, by slight abuse of language, itself a tubular neighborhood of $S$ in $M$.

Following Mather \cite{MatNTS} we define the \emph{projection} $\pi_S: T_S \to S$ as the composition $\pi_E \circ \varphi^{-1}$
and the \emph{tubular function} $\varrho_S : T_S \to [0,\infty)$ as the function which maps a point $x \in T_S$ to
$\big\langle \varphi^{-1} (x),  \varphi^{-1} (x) \big\rangle$. 

If $N$ is a a second $G$-manifold and $f:M \to N$ a $G$-equivariant smooth map, then a tubular neighborhood 
$\sfT=(E,\varepsilon,\varphi)$ is said to be \emph{compatible} with $f$ if $f \circ \pi_S = f|_{T_S}$.

\begin{example}
  Let $\eta$ be a $G$-invariant riemannian metric on $M$ and $\pi^N: N\to S$ the normal bundle 
  of $S$ in $M$. Identify $N$ with the orthogonal complement of $TS$ in $T_SM$ via $\eta$. 
  Then there exists an open neighborhood $U$ of the zero section of $N$ such that 
  the exponential function $\exp : U \to M$ is a $G$-equivariant 
  diffeomeorphism onto an open neighborhood of $S$ in $M$.
  Since $G$ is compact there exists a $G$-invariant continuous function 
  $\varepsilon $  such 
  $B(\varepsilon,N) := \{ v \in E \mid \langle v,v \rangle < \varepsilon (\pi_N(v)) \} \subset U$
  The inner product $\langle -,- \rangle$ on $N$ hereby is the one induced by the riemannian metric 
  $\eta$. The triple $\sfT = (N, \varepsilon, \exp|_{B(\varepsilon,N)})$ now is a 
  $G$-equivariant tubular neighborhood of $S$ in $M$.
\end{example}

\begin{example}
  Consider the action of the Lie group $ \SU(2)$ on itself by conjugation.
  Since a matrix $g \in \SU(2)$ is diagonalizable it is conjugate to a matrix
  of the form
  $\begin{pmatrix} e^{i\lambda} & 0 \\ 0 & e^{-i\lambda}\end{pmatrix}$
  with $\lambda \in [0,\pi]$ uniquely determined. The conjugacy classes of
  $ \SU(2)$ can therefore be labelled by the elements of the interval $[0,\pi]$.
  Denote the corresponding map $ \SU(2) \to [0,\pi]$ by $\lambda$.
  Note that $\lambda$ is continuous but not smooth. 
  Next observe that the \emph{latitude} map
  $\ell = \frac12  \tr : \SU (2) \to \R$ is smooth, 
  equivariant and has image $[-1,1]$. Moreover, if $g$ has eigenvalues
  $e^{\pm i\lambda}$, then $\ell (g) =\cos \lambda$, which in other words means
  that the level map labels the conjugacy classes of $\SU(2)$ as well.
  Now equip   
  $M := \SU(2)$ with the stratification by orbit types. There are three of
  those, namely
  \begin{align} \nonumber
    M_{(\SU(2))} & = \{ \pm 1 \} = \ell^{-1} \big(\{\pm 1 \}) = \lambda^{-1} ( \{ 0,\pi \} \big), \\ \nonumber
    M_{(\mathsf{O}(2))} & =  \ell^{-1} ( \{ 0 \} )  = \lambda^{-1} ( \{ \pi/2 \} ), \text{ and }  \\ \nonumber
    M_{(\mathsf{U} (1))} & =
    \ell^{-1} \big( (-1, 1) \setminus \{ 0 \} \big) = \lambda^{-1} \big( (0,\pi)\setminus \{\pi/2\} \big) \ . 
  \end{align}
  When endowing $[-1,1]$ with the stratification given by the open subset
  $(-1, 1) \setminus \{ 0 \}$ and the
  discrete subset $\{0,\pm 1\}$ the level map $\ell$ becomes a smooth
  stratified equivariant submersion. 

  Let us now explicitly describe a system of equivariant control data on $\SU(2)$. Put
  \begin{align} \nonumber
    T_{(\SU(2))} & = \lambda^{-1} \big( [0,\pi/4) \cup ( 3 \pi/4 , \pi/4] \} \big), \\ \nonumber
    T_{(\mathsf{O}(2))} & = \lambda^{-1} \big( ( \pi/4 ,3 \pi/4 ) \big), \text{ and }  \\ \nonumber
    T_{(\mathsf{U} (1))} & = M_{(\mathsf{U} (1))} = \lambda^{-1} \big( (0,\pi)\setminus \{\pi/2\} \big) \ .
  \end{align}
  Then $T_{(\SU(2))}, T_{(\mathsf{O}(2))}$ and $T_{(\mathsf{U} (1))}$  are tubular neighborhoods of
  $M_{(\SU(2))},  M_{(\mathsf{O}(2))}$ and $M_{(\mathsf{U} (1))}$, respectively. 
  The projections and tubular functions are defined for diagonal
  $g = \begin{pmatrix} e^{i\lambda} & 0 \\ 0 & e^{-i\lambda}\end{pmatrix} \in \SU(2)$
  as follws and then extended equivariantly:
  \begin{align}\nonumber
    \pi_{(\SU(2))} (g) & =
       \begin{cases}
          \begin{pmatrix} 1 & 0 \\0 & 1 \end{pmatrix}& \text{if } \lambda \in [0,\pi/4) \ ,\\ 
          \begin{pmatrix} -1 & 0 \\0 & -1 \end{pmatrix}& \text{if } \lambda \in (3\pi/4,\pi] \ ,
       \end{cases}  \\ \nonumber
    \pi_{(\mathsf{O}(2))} (g)  & = \begin{pmatrix} i & 0 \\0 & -i \end{pmatrix}  \quad
       \text{if } \lambda \in ( \pi/4 ,3 \pi/4 ) \ ,  \\ \nonumber 
  %\end{align}
  %\begin{align}\nonumber
    \varrho_{(\SU(2))} (g) & = \begin{cases}
          \lambda^2 & \text{if } \lambda \in [0,\pi/4) \ , \\ 
          (\lambda - \pi)^2 & \text{if } \lambda \in (3\pi/4,\pi] \ ,
       \end{cases}  \\ \nonumber
    \varrho_{(\mathsf{O}(2))} (g) & = (\lambda - \pi/2)^2\quad \text{if } \lambda \in ( \pi/4 ,3 \pi/4 ) \ .
  \end{align}
  The projection $\pi_{(\mathsf{U} (1))}$ and tubular function $\varrho_{(\mathsf{U} (1))}$ are the identity
  and zero map, respectively. One checks that the thus defined data give a
  system of equivariant control data for $\SU(2)$ carrying the conjugate action.

  In Theorem \ref{thm:control-data-existence} we will establish a general scheme to construct equivariant
  systems of control data under the assumption that the underlying stratified space is (B) regular and that
  the stratification is compatible with the group action. 
\end{example}

Given two tubular neighborhoods $\sfT = (E, \varepsilon, \varphi)$  and $\sfT' = (E', \varepsilon', \varphi')$ 
an \emph{isomorphism} between $\sfT$ cand $\sfT'$ consists of an isometric vector bundle isomorphism 
$\psi : E \to E'$ and a continuous map $\delta : S \to (0,\infty)$  such that 
$\delta \leq \min (\varepsilon , \varepsilon')$ and such that 
$ \varphi' \circ \psi|_{B(\delta,E)} = \varphi|_{B(\delta,E)}$.  We denote such an isomorphism briefly by 
$\psi: \sfT \sim \sfT'$. 

In addition to isomorphisms of equivariant tubular neighborhoods there is a corresponding equivariant version of diffeotopies
on a $G$-manifold $M$. By a $G$-equivariant diffeotopy on $M$ we understand a smooth map $H :  M \times I \to M$, where 
$I$ denotes the interval $[0,1]$, such that each of the maps $H_t :M \to M$, $x \mapsto H(x,t)$ with $t\in I$ 
is $G$-equivariant diffeomorphism and such that $H_0 = \id_M$. 
If in addition a $G$-equivariant map $f: M \to N$ is given, then a $G$-equivariant diffeotopy $ H :  M \times I \to M$ is said 
to be \emph{compatible} with $f$ if $f (H(x,t)) = f(x)$ for all $x\in M$ and $t\in I$. 
For every diffeotopy $H : M \times I \to M$ one calls the set 
$ \supp H := \overline{\{ x \in M \mid \exists t\in I : H(x,t) \neq x \}}$ the \emph{support} of $H$. 

If $h: (M,S) \to (M',S')$ is a $G$-equivariant diffeomorphism between pairs of $G$-manifolds and $G$-submanifolds
one defines the push-forward tubular neighborhood $h_* \sfT$ of $S'$ in $M'$ of a tubular neighborhood $\sfT = (E,\varepsilon,\varphi)$ 
of $S$ in $M$  by  $h_* \sfT = \big( (h^{-1})^*E, \varepsilon \circ h^{-1}, h \circ \varphi \big)$.

\subsection{The Equivariant Submersion Theorem}
Locally in charts, every submersion looks like a linear projection. The following equivariant versions of this result appear 
to be folklore.  
\begin{proposition}
\label{prop:LocalFormEquivariantSubmersion}
  Let $f:M \to P$ be a $G$-equivariant submersion between $G$-manifolds $M$ and $P$.
  Let $x\in M$ be a point, $H=G_x$ be the isotropy group of $x$ and $K=G_{f(x)} \supset H$ the one of 
  $f(x)$. Then there exist a finite dimensional orthogonal  $K$-representation space $N$, a finite dimensional orthogonal  $H$-representation
  space $W$,  a $K$-invariant open convex neighborhood of the origin $B\subset N$,
  an $H$-invariant open convex neighborhood of the origin $C\subset W$,
  $G$-equivariant open embeddings $\Theta: G\times_H (B \times C) \hookrightarrow M$ and $\Psi: G\times_K B \hookrightarrow P$ such that 
  $\Theta \big( [e,0]_{G\times_H (B \times C)} \big) = x$,  $\Psi \big( [e,0]_{G\times_K B}  \big) = f(x)$ and such that the following 
  diagram commutes, where $\pi : B \times C \to C$ is projection onto the first factor and 
  $\overline{\id_G \times \pi}$ maps $[g,(v,w)]_{G \times_H (B\times C)}$ to $[g,v]_{G\times_H B}$. 
  
  \begin{equation}
  \label{Dia:LocalFormEquivariantSubmersion}
  \xymatrix{
  G\times_H (B \times C) \ar[rr]^{\overline{\id_G\times \pi}} \ar[d]_\Theta & \hspace{5mm}& G\times_KB \ar[d]^\Psi \\
  M \ar[rr]_f & \hspace{5mm}& P
  }
  \end{equation}
 Moreover, the dimensions of the manifolds and representation spaces fulfill the relations 
 $\dim P = \dim G - \dim K + \dim N$ and $\dim M = \dim G - \dim H + \dim N + \dim W$.  
\end{proposition}
\begin{proof}
  First choose  a $G$-riemannian metric $\varrho$ on $P$, identify for every $p \in P$ the normal space 
  $N_{p} := T_{p}P/T_{p} Gp$ with the orthogonal complement of $T_{p}Gp$ in $T_{p}P$, put $N:=N_{f(x)}$ and let 
  $B \subset N$ be an 
  open ball around the origin of radius smaller than the injectivity radius of the exponential function $\exp^\varrho $ 
  (with respect to $\varrho$) over the orbit $Gf(x)$. Then, $N$ is an orthogonal $K$-representation space and,
  by the Slice Theorem \cite[p.~139]{KosCGTL}, \cite{PalESANCLG}, \& \cite[Sec.~II.4]{BreICTG},
  the subset $Z =\exp^\varrho (B) \subset P$ is a $K$-invariant submanifold and the map 
  \[ 
    \Psi: G\times_K B \to P, \: [g,v]_{G\times_K D} \mapsto g \exp^\varrho v
  \]
  a $G$-equivariant diffeomorphism onto an open  neighborhood of $Gf(x)$ in $P$.
  Following \cite[Sec.~II.4]{BreICTG}, we call such a $Z\subset P$ a \emph{slice} 
  through $f(x)$.  
  Next choose a $G$-equivariant bundle $\calH\to M$ complementary to the vertical bundle $\calV=\ker Tf \to M$ in $TM$ 
  and call it the horizontal bundle. For each $v \in B$ let $\widetilde{\gamma}_v : [0,1] \to M$ be the horizontal 
  lift of the geodesic $\gamma_v : [0,1] \to P$, $t\mapsto \exp^\varrho (tv)$ such that $\widetilde{\gamma}_v (0) = x$.
  Let $\widetilde{Z} = \{ \widetilde{\gamma}_v (1) \mid v \in B \}$. Then   
  $\widetilde{Z}$ is a submanifold of $M$ with a 
  chart given by the inverse of $B \to \widetilde{Z}$, $v \mapsto  \widetilde{\gamma}_v (1)$. Moreover,  $\widetilde{Z}$ 
  is $H$-invariant since the action by $h \in H$ maps the geodesic $\gamma_v$ for $v \in B$  to the geodesic 
  $\gamma_{hv}$ and the lift of $\gamma_v$ to the lift of $\gamma_{hv}$. 
  After possibly shrinking $B$ and with it $Z$ and $\widetilde{Z}$ the map $G\times_H \widetilde{Z} \rightarrow M$ becomes 
  a $G$-equivariant embedding and the composition  $G\times_H \widetilde{Z} \rightarrow M \overset{f}{\rightarrow} P$ 
  a submersion. Now consider the fiber $F \subset M$ of $f$ through the point $x$. Then $F\subset M$ is a submanifold 
  which is invariant under the action of $K$. Choose a $K$-invariant riemannian metric $\eta$ on $F$ and let 
  $W := T_xF/T_xKx$ be the normal space at $x$ and identify it with the orthogonal complement of $T_xKx$ in $T_xF$. 
  Observe that $W$ is an orthogonal $H$-representation space since $H = K_x$. For a suficiently small open ball 
  $C \subset W$ around the origin the set $Y := \exp^\eta (C)$ is a slice of the $K$-manifold $F$ through $x$.
  Now let $\widetilde{\gamma}_{y,v}$ for each $y\in Y$ and $v \in B$ be the horizontal lift of $\gamma_v$ such that 
  $\widetilde{\gamma}_{y,v}(0)=y$. 
  Then, after possibly shrinking $B$ and $C$ (together with the corresponding slices) the map  
  $\theta : B \times C \hookrightarrow M$, $(v,w)\mapsto \widetilde{\gamma}_{\exp^\eta w,v}(1) $ is
  an embedding since $T_0\theta$ is the linear injection 
  $N \times W \hookrightarrow \calH_x \times \calV_x\cong T_xM$ which maps $(v,w)$ to the unique pair 
  $(v^\textup{h},w) \in  \calH_x \times \calV_x$ such that $T_x(v^\textup{h})=v$. 
  Moreover, the embedding $\theta$ is $H$-equivariant since the horizontal lift is $H$-equivariant. 
  Hence the set $X= \{  \widetilde{\gamma}_{y,v}(1)\mid v \in B, \: y \in Y \}$
  is an $H$-invariant submanifold of $M$ transversal to the orbit $Gx$. 

  Consider now the $G$-manifold $G \times_H (B \times C)$ and define
  \[
    \Theta : G \times_H (B \times C) \to M, \: [g, (v,w)]_{G \times_H (B \times C)} \mapsto  g \theta (v,w) \ .
  \]
  This map is well-defined by equivariance of $\theta$. The restriction of $\Theta$ to the zero section of  
  $G \times_H (B \times C)$, which is canonically diffeomorphic to $G/H$, is a diffeomeorphism
  onto to the orbit $Gx$ since $H=G_x$. Moreover, the image of $\{ e \} \times (B \times C)$ under $\Theta$
  is the manifold $X$, so the image of $\Theta$ equals $GX$. Now recall that the orbit $Gx$ and $X$ 
  are transverse. By $G$-equivariance of $\Theta$ it follows that $T_{gx}\Theta$ is an isomorphism for all $g\in G$. 
  After possibly shrinking $B$ and $C$ the map $\Theta$ therefore is a diffeomorphism onto an open neighborhood of 
  the orbit $Gx$ by the Implicit Function Theorem.  In other words, $X$ is a slice of $M$ at $x$.      

  Finally  let $\psi : B \to P$ be the embedding $v \mapsto \exp^\varrho v$.
  The image of $\psi$ then is the slice $Z$, and $ \psi^{-1} \circ f \circ \theta = \pi$ by definition of 
  $\theta$ via horizontal lifts. This entails that the diagram commutes and the proposition is proved. The dimension relation 
  follows from the definition of $N$ and $W$ and the Slice Theorem. 
\end{proof}

\begin{corollary}
\label{cor:LocalFormEquivariantSubmersionSubmanifold}
  Let $f:M \to P$ be a $G$-equivariant map between $G$-manifolds $M$ and $P$, and let $S\subset M$ be a
  $G$-invariant submanifold such that the restriction $f|_{S}:S \to P$ is a submersion. 
  Let $x\in S$ be a point, $H=G_x$ be the isotropy group of $x$ and $K=G_{f(x)} \supset H$ the one of 
  $f(x)$. Then there exist an orthogonal  $K$-representation space $N$,   
  orthogonal  $H$-representation spaces $W$ and $W'$, a $K$-invariant open convex neighborhood of the 
  origin $B\subset N$,  $H$-invariant open convex neighborhoods of the origin $C\subset W$  and $D \subset W'$,
  and finally $G$-equivariant open embeddings $\Phi : G\times_H (B \times C \times D) \hookrightarrow M$ and 
  $\Psi: G\times_K B \hookrightarrow P$ such that 
  $\Phi \big( [e,0]_{G\times_H (B\times C \times D)} \big) = x$,  $\Phi \big( [e,0]_{G\times_K B}  \big) = f(x)$ and such 
  that the following two properties hold true.
  \begin{enumerate}[\rm (1)]
  \item \label{ite:ImageInSubmanifold}
  The map 
  \[
    \Theta : G \times_K (B\times C) \to M, \: [g,(v,w)]_{G \times_K (B\times C)} \to \Psi \big( [v,w,0]_{G \times_K (B\times C \times D)} \big)
  \]
  has image in $S$ and comprises a $G$-equivariant diffeomorphism onto an open neighborhood of the 
  orbit $Gx$ in $S$. 
  \item  
  With $\Pi : B \times C \times D \to B$ denoting projection onto the first factor the diagram below commutes. 
  \begin{equation}
  \label{Dia:LocalFormEquivariantSubmersionSubmanifold}
  \xymatrix{
  G\times_H (B \times C \times D) \ar[rr]^{\overline{\id_G\times \Pi}} \ar[d]_\Phi & \hspace{5mm}& G\times_KB \ar[d]^\Psi \\
  M \ar[rr]_f & \hspace{5mm}& P
  }
  \end{equation}
  \end{enumerate}
\end{corollary}

\begin{proof}
 Choose a $G$-invariant riemannian metric $\varrho$ on $P$ and a $G$-equivariant horizontal bundle 
 $\calH \to P$ such that the restriction $\calH|_S$ is a $G$-equivariant subbundle of $TS$. 
 Let $N = T_{f(x)}P/T_{f(x)}Gf(x)$ and choose $B \subset N$ as in the proof of the proposition.
 Then $Z=\exp^\varrho B$ is a slice and $\Psi :G \times_K B\to B$ defined as above a $G$-equivariant diffeomorphism.    
 Now let $F\subset M$ be the fiber of the submersion $f$ through $x$. Then $F$ is a $K$-invariant submanifold of $M$ and
 $F \cap S$ a $K$-invariant submanifold of $S$.  Choose a $K$-invariant riemannian metric $\eta$ on $F$.  
 Put $W =T_x(F\cap S)/T_xKx$ and identify it with the orthogonal complement of $T_xKx$ in $T_x(F\cap S)$. 
 Let $W'$ be the orthogonal complement of $T_x(K\cap S)$ in $T_xF$ and choose small enough open convex neighborhoods of the 
 origin $C \subset W$ and $D\subset W'$. Then $Y = \exp^\eta (C\times D)$ is a slice through $x$ of the $K$-manifold $F$. 
 Let $\gamma_v$ for $v\in B$ the geodesics as in the proof of the proposition and denote for $y\in Y$ and $v\in B$ by
 $\widetilde{\gamma}_{y,v}$ the horizontal lift of $\gamma_v$ in $M$ which starts at $y$. 
 Now define $\varphi : B \times C \times D \to M$ by $(v,w,z)\mapsto \widetilde{\gamma}_{\exp^\eta(w,z),v}$. 
 Then 
 \[
  \Phi :  G\times_H (B \times C \times D) \hookrightarrow M , \: [g, (v,w,z)]_{ G\times_H (B \times C \times D)} \mapsto
  g \varphi (v,w,z) 
 \] 
 is a $G$-invariant diffeomorphism onto an open neighborhood of $Gx$ in $M$ and the diagram 
 \eqref{Dia:LocalFormEquivariantSubmersionSubmanifold} commutes by the proof of the proposition.
 By assumptions on the horizontal bundle $\calH$ and the construction of $\Theta$ property 
 \eqref{ite:ImageInSubmanifold} holds true. 
\end{proof}

\subsection{Uniqueness and existence of equivariant tubular neighborhoods}

For the construction of $G$-equivariant control data one needs stronger 
versions of existence and uniqueness results of $G$-equivariant tubular neighborhoods.  
In the following we prove equivariant versions of \cite[Prop.~6.1]{MatNTS} and 
\cite[Prop.~6.2]{MatNTS}.

\begin{theorem}[Uniqueness of equivariant tubular neighborhoods]\label{thm:tubular-uniqueness}
\label{thm:uniqueness-equivariant-tubular-neighborhoods}
  Let $M$, $P$ be smooth $G$-manifolds, $S \subset M$ a closed $G$-invariant smooth submanifold, and 
  $f: M\to P$ a $G$-equivariant smooth map such that the restriction $f|_{S} :S \to P$ is a submersion. 
  Assume that $\sfT_0$ and $\sfT_1$ are two $G$-equivariant tubular neighborhoods of $S$ in $M$ and 
  that they are compatible with $f$. 
  Further assume that $U\subset S$ is a $G$-invariant relatively open subset and that 
  $\psi_0 :\sfT_0|_U \to \sfT_1|_U$ is an  isomorphism of $G$-equivariant tubular 
  neighborhhoods over $U$. 
  Let $A,Z \subset S$ be two $G$-invariant relatively closed subsets such that $A \subset U$ and let 
  $V\subset M$ be a $G$-invariant open neighborhood of $Z$ in $M$.
  Then  there exists a $G$-equivariant diffeotopy $H: M \times I \to M$ which leaves $S$ invariant, is compatible 
  with $f$ and has support in $V\setminus A$ such that  the tubular neighborhoods $(H_1)_* \big(\sfT_0|_{A\cup Z}\big)$ 
  and $\sfT_1|_{A\cup Z}$ are isomorphic.
  If $O\subset M \times M$ is a $G$-invariant neighborhood of the diagonal, one can choose $H$ such 
  that $(H_t(x),x) \in O$ for all $t \in I$ and $x\in M$. 
  Finally, the isomorphism $\psi : (H_1)_* \big(\sfT_0|_{A\cup Z}\big) \to \sfT_1|_{A\cup Z}$ is $G$-equivariant 
  and can be constructed so that $\psi|_A = \psi_0|_A$.  
\end{theorem}

\begin{proof}
  Our proof adapts Mather's argument in \cite[Proof of Prop.~6.1]{MatNTS} to the $G$-equivariant case. 

  {\it Step 1.} We first consider the local $G$-equivariant case as stated in Corollary \ref{cor:LocalFormEquivariantSubmersionSubmanifold}.
  So we assume for now the following:
  \begin{enumerate}[(1)]
  \item  $P$ is of the form $G \times_{G_{f(s)}} B$, where $G_{f(s)} \subset G$ is the isotropy group of some point 
         $f(s)$ with $s\in S$ and $B$ is an open convex neighborhood of the origin of some euclidean space
          $\R^k$ carrying an orthogonal $G_{f(s)}$-representation. The point $f(s)$ is then identified with $[e,0] \in G \times_{G_{f(s)}} B$.
  \item  $S$ is equivariantly diffeomorphic to an associated bundle of the form $G \times_{G_s} (B \times C)$, 
         where $G_s \subset G$ is the isotropy group of  $s\in S$ 
         and $C$ is an open convex neighborhood of the origin of some euclidean space $\R^l$ carrying an orthogonal 
         $G_s$-representation. Under the corresponding diffeomorphism  the point $s$ can be identified with 
         $[e,0]\in G \times_{G_s} (B\times C)$. Note that $G_s \subset G_{f(s)}$.
         % and the fiber over any point  $x \in S$ with the manifold $B$. 
  \item  $M$ is equivariantly diffeomorphic to an associated bundle of the form $G \times_{G_s} (B \times C \times D)$,
         where $D$ is an open convex neighborhood of the origin of a euclidean space $\R^m$ with an orthogonal $G_s$-representation
         and where the $G_s$-action on $B \times C \times D$ is the diagonal action. 
         %The fiber over $s$ in $M$ is identified with $B \times C$. 
  \item  Under these identifications $f: M \to P$ coincides with the $G$-equivariant map 
         $G \times_{G_s} (B \times C \times D) \to G\times_{G_{f(s)}}B$ which maps
         $[g,(v,w,z)]$ to $[g,v]$. 
         So for every $x = [g,(v,w,z)] \in M$ the fiber through  $x$ in $M$ coincides
         with $F_x := \big[ \{ g \} \times \{ v\} \times C \times D\big]$, the image of 
         $\{ g \} \times \{ v\} \times C \times D$ in  $G \times_{G_s} (B \times C \times D)$.
  \end{enumerate}
  In addition to this we also assume for the moment that $Z$ is compact.
  
  Since $G$ is compact, there exists a bi-invariant riemannian metric $\mu$ on $G$. The spaces $B,C,D$ all carry 
  natural invariant metrics induced by the ambient euclidean spaces. Denote by $\eta$ and $\varrho$ the
  induced $G$-invariant riemannian metrics on $M$ and $P$, respectively.  With these metrics,
  $f$ then becomes a riemannian submersion. Actually, the fibers of this riemannian submersion are even totally geodesic
  by construction of $\eta$. Now assume that $x$ and $y$ are points of $M$ which are both in the same fiber $F_x$. 
  Then $x = [g,(v,w,z)]$ and  $y = [g,(v,w',z')]$ for some $g\in G$, $v\in B$, $w,w'\in C$ and $z,z'\in D$. 
  The unique geodesic connecting $x$ with $y$ then is given by 
  \[
   \gamma_{x,y}(t) =  [g,(v,(1-t)w+t w',(1-t)z+tz')] \quad \text{for all $t\in I$, where $I=[0,1]$}.  
  \]
  Note for later that $\gamma_{x,y}$ completely runs within the fiber $F_x$. 

  % By construction, the orthogonal complement of the tangent bundle $TS$ in $T_SM := \bigcup_{x\in S} T_xM$ 
  % then is contained in the kernel bundle $\ker T_S f := \bigcup_{x\in S} \ker T_xf$.    
  Denote by $\pi_N : N \to S$ the normal bundle of $S$ in $M$ that is $N_x = T_xM/T_xS \cong T_0D \cong \R^m$ for all $x\in S$.
  Via the riemannian metric $\eta$ one can identify $N$ with the subbundle of $T_SM$ orthogonal to $TS$.
  By assumptions and construction of the riemannian metric $\eta$ one has $N \subset \ker T_Sf$. 

  Now observe that for $i=0,1$ the map
  \[
     \alpha_i : E_i \to N, \quad v \mapsto T\varphi_i (v) + T_{\pi_{E_i} (v)} S
  \]
  is a vector bundle isomorphism. Hereby we have identified $E_i$ with the vertical subbundle of $TE_i$ restricted to the zero section. 
  By assumptions $\alpha_i$ is an isomorphism of $G$-bundles, hence $\alpha := \alpha_1^{-1} \circ \alpha_0 : E_0 \to E_1$ is one, too.
  Note that for $x \in U$,  $\alpha_x  : E_{0,x} \to E_{1,x}$ coincides with  $\psi_{0,x}  : E_{0,x} \to E_{1,x}$.
  By uniqueness of the polar decomposition there exists a unique $G$-equivariant vector bundle automorphism 
  $\beta : E_1 \to E_1$ such that for every $x\in S$ the linear map $\beta_x: E_{1,x} \to E_{1,x}$ is positive definite  and 
  $\psi_x := \beta_x \circ \alpha_x : E_{0,x} \to E_{1,x}$ an orthogonal transformation. 
  Then
  \[
    \xi_t := (1-t)\alpha + t \psi : E_0 \to E_1    
  \] 
  is an isomorphism for every $t\in I$ which over $U$ coincides with $\psi_0$. 
  After possibly lessening $\varepsilon_1$ and $\varepsilon_0$, where both stay $G$-invariant and positive, 
  the set $T := T_{0,S} \cap T_{1,S}$ is a $G$-invariant open neighborhood of $S$ in $M$   over which the maps 
  \[
    q_t : \: T \to M, \: x \mapsto \varphi_1 \circ   \xi_t \circ \varphi_0^{-1} (x) 
  \]
  are well-defined and open $G$-equivariant embeddings for every $t\in I$. Note that over $S$ each $q_t$ is the identical embedding
  and that each $q_t$ acts as identity over some open $G$-invariant neighborhood $U' \subset T$ of $A$. 
  Moreover, each $q_t$ is compatible with $f$ since both $\sfT_0$ and $\sfT_1$ are compatible with $f$. 
  Put $V_1 = T \cap V$ and observe that $Z \subset V_1$. By compactness of $Z$ there exists a $G$-invariant
  open neighborhood $V_2$ of $Z$ which is relatively compact in $V_1$ and which satisfies $V_2 \subset q_t (V_1)$ 
  for all $t\in I$. Next choose a smooth $G$-invariant function $\chi : M \to I$ with compact support in $V_2$ 
  such that $\chi$ is identically $1$ over a $G$-invariant neighborhood of $Z$ in $V_2$. 
  Define $Q_{s,t}: M \to M$ for $s,t\in I$ by 
  \[
   Q_{s,t} (x) = 
   \begin{cases}
     \gamma_{x,q_t \circ q_s^{-1} (x)} \big( \chi(x) \big) & \text{if $x\in V_2$}, \\
     x & \text{if $x\in M \setminus V_2$},
   \end{cases}
  \]
  Since the $q_t$ are compatible with $f$, the geodesic $\gamma_{x,q_t \circ q_s^{-1} (x)}$ is well-defined for every $x\in V_2$. 
  Hence the $Q_{s,t}$ are well-defined as well and also compatible with $f$. Next observe that by construction $Q_{t,t}$ is the 
  identity map for all $t \in I$ and that there is a compact $G$-invariant subset containing the support of $Q_{s,t}$ for all $s,t\in I$. 
  Hence there is some $\delta > 0$ such that $Q_{s,t}$ is a diffeomorphism for all $s,t$ with $|s-t|< \delta$. 
  From here on we can follow Mather's treatment of the local case in \cite[Prop.~6.1]{MatNTS}  almost literally.
  Choose a positive integer $n$ such that $\frac 1n < \delta$. Put
  \[
  \widetilde{H}_t = Q_{0,\frac t n} \circ  Q_{\frac tn, \frac{2t}{n}} \circ  \ldots   \circ  Q_{\frac{(n-1)t}{n},t} \ .
  \]
  Then $\widetilde{H}$ is a $G$-equivariant diffeotopy, compatible with $f$, and leaves $S$ fixed. Since the $q_t$ acts as identity over $U'$,
  $\widetilde{H}_t$ does so, too.  Moreover, $\widetilde{H}$ coincides by construction with $q_1 \circ q_0^{-1}$ over some 
  sufficiently small $G$-invariant open neighborhood of $Z$ in $V_2$. Hence  $\widetilde{H}$ coincides with $q_1 \circ q_0^{-1}$ over
  $U'\cup V'$. Furthermore $\widetilde{H}\circ q_0\circ \varphi_0 = q_1 \circ \varphi_0 = \varphi_1 \circ \psi $
  over the $G$-invariant neighborhood $\varphi_0^{-1} (U'\cup V')$ of $A \cup Z$ in $E_0$. 
  Therefore, $\psi$ is an isomorphism between $(\widetilde{H}_1 q_0)_* \sfT_0|_{A\cup Z} $ and $\sfT_1|_{A\cup Z}$. 
  Moreover, the support of $\widetilde{H}$ is contained in $V_2 \subset V$ by construction. 
  Finally, by requiring that the cut-off function $\chi$ has support in a sufficiently small $G$-invariant open neighborhood of $Z$ 
  one can achieve that with regard to the compact-open topology $\widetilde{H}_t$ is uniformly in $t\in I$ as close to the identity
  map as one wishes. 
  
  By the following step there exists, after possibly shrinking $U'$ and $V'$, a $G$-equivariant diffeotopy $\widehat{H}$ of $M$ which is compatible 
  with $f$ and leaves $S$ invariant such that for all $t\in I$  the diffeomorphisms $\widehat{H}_t$ act as identity over $U'$ 
  and such that $\widehat{H}_t =q_0$ over $V'$.  The map $H : M \times I \to M$, $(x,t) \mapsto \widetilde{H}_t \circ \widehat{H}_t (x)$ then 
  is a $G$-equivariant diffeotopy with all the required properties.

  {\it Step 2.} Here we show that there exists a $G$-equivariant diffeotopy $\widehat{H}$ of $M$ with compact support which is compatible with $f$, 
  leaves $S$ invariant, acts as identity over a sufficiently small open $G$-invariant neighborhood $U'$ of $A$ and coincides over a sufficiently 
  small open $G$-invariant neighborhood $V'$ of $Z$ with $q_0$. Note that for  the non-equivariant case the existence of such a diffeotopy 
  $\widehat{H}$ has been claimed in the proof of \cite[Prop.~6.1]{MatNTS} with the argument left to the reader. Since the equivariant case is more 
  subtle, we present a proof here which obviously covers Mather's claim, too. 

  Observe that for all $x \in T$ the image $q_0(x)$ lies in the fiber $F_x$ by construction. 
  Hence the geodesic $\gamma_{x,q_0(x)}$ is well-defined and fully 
  runs in $F_x$. Now put $K(x,t) = \gamma_{x,q_0(x)} (t)$ for all $x\in T$ and $t\in I$. After possibly shrinking $T$, 
  $\widehat{K}: T \times I \rightarrow M \times I$, $(x,t)\mapsto \big( K(x,t),t\big)$  is an open embedding since $q_0$ acts as identity over $S$,
  one has $T_xq_0 = \id_{T_xM}$ for all $x \in S$ and finally since $I$ is compact. Moreover, $K_t$ acts as identity for all $x\in U'$
  where $U'$ is a $G$-invariant open neighborhood of $A$ in $T$ over which $q_0$ acts as identity. Finally, $K$ is compatible with $f$
  since $q_0$ is.   
  Now define the time-dependent vector field $X_K: T \times I \to TM$ by 
  \[
   X_K(x,t) = \left.\frac{\partial}{\partial s}\right|_{s=0} K(x,t +s) \ .
  \]
  Note that over $U' \times I $ the vector field $X_K$ vanishes and that $X_K$ is $G$-equivariant.
  Next choose a sufficiently small relatively compact $G$-invariant open neighborhood $V'$ of $Z$ and a 
  non-negative $G$-invariant smooth function $\delta : M \to [0,1]$ which is identical to $1$ over $V'$ and 
  has compact support in $T$.   
  Define the time-dependent vector field $X:M \times I \to TM$ by 
  \[ 
   (x,t)\mapsto 
   \begin{cases}
     \delta(x) \, X_K(x,t) & \text{for $x \in T$}, \\
     0 & \text{for $x \in M \setminus T$}. \\
   \end{cases}
  \]
  Then $X$ is a time-dependent $G$-equivariant vector field on $M$ with compact support. By \cite[Chap.~8, Thm.~1.1]{HirDT} 
  it generates a diffeotopy $\widehat{H}:M \times I \to M$. The diffeotopy is $G$-equivariant since $X$ is, and is compatible with $f$ since
  $X$ is tangent to the fibers of $f$ by construction. Since $\delta$ has compact support, $\widehat{H}$ has so too. 
  Over $V'$ the diffeotopy $\widehat{H}$ coincides with $K$, hence one obtains in particular that $\widehat{H}_1|_{V'} = q_0|_{V'}$.
  Over $U'$, each $\widehat{H}_t$ acts as identity for every $t\in I$. This finishes Step 2.
 
  {\it Step 3.} Let us pass to the general case, now.  Here we follow closely  \cite[Prop.~6.1]{MatNTS}.
  By the Equivariant Submersion Theorem and Corollary \ref{cor:LocalFormEquivariantSubmersionSubmanifold} 
  there exists for every $x\in S$ an open relatively compact $G$-invariant
  open neighborhood $W_x$ of $x$ in $M$ together with $G$-equivariant open embeddings called  
  \emph{equivariant charts} $\Phi_x :  W_x \hookrightarrow G \times_{G_x} \R^{p+k+l} $  
  and $\Psi_x : f(W_x) \hookrightarrow   G\times_{G_{f(x)}} \R^p$, 
  where $\R^p$ carries an orthogonal $G_{f(x)}$-representation and $\R^k$  and $\R^l$ orthogonal $G_x$-representations,
  such that the following conditions hold true:
  \begin{enumerate}[(1)]
  \item \label{cond:one}
   The image of $\Phi_x$ is of the form $G\times_{G_x} (B \times C \times D)$ with 
   $B \subset \R^p$, $C \subset \R^k$, and $D\subset \R^l$ open convex neighborhoods of the origin, and 
   $\Psi_x \big( f(W_x) \big) = G \times_{G_{f(x)}} B$. 
  \item 
   One has 
   $W_x \cap S = \Phi_x^{-1} \big( G \times_{G_x} (B \times C \times \{ 0 \})\big) = \Phi_x^{-1}\big( G \times_{G_x} \R^{p+k} \big)$.     
  \item \label{cond:three}
   The diagram
   \begin{displaymath}
   \xymatrix{
     W_x \ar[rr]^{\Phi_x} \ar[d]_f & \hspace{5mm} & G\times_{G_x} \R^{p+k+l} \ar[d]^{\overline{\id_G \times \Pi}}\\
     f(W_x) \ar[rr]_{\Psi_x} & \hspace{5mm} & G\times_{G_{f(x)}} \R^p
   }
  \end{displaymath}
   commutes, where $\Pi$ is projection onto the first $p$ coordinates.
  \end{enumerate}
  After possibly shrinking the $W_x$ one can achieve that 
  \begin{equation}
    \label{eq:condition}
    \begin{split}
      W_x \cap A \neq \emptyset & \implies  W_x \subset V \\
      W_x \cap Z \neq \emptyset & \implies  W_x \cap S \subset U \ . 
    \end{split}
  \end{equation}
  The family ${M \setminus S} \cup \{ W_x \}_{x \in S}$ then is covering of $M$ by $G$-invariant open subsets. Since the orbit space $M/G$ is separable and 
  paracompact one can find a locally finite countable refinement  ${M \setminus S} \cup \{ W_i \}_{i \in \NBbb^*}$ with each $W_i$ being $G$-invariant, open in $M$ and 
  contained in $W_{x_i}$ for some $x_i \in S$. Moreover, the $W_i$ are so that there exist equivariant charts 
  $\Phi_i : W_i \to G \times_{G_i} \R^{p+k+l}$ and 
  $\Psi_i: f(W_i) \to  G \times_{G_{f(i)}} \R^p$ fulfilling conditions 
  \eqref{cond:one} to \eqref{cond:three}.
  Following Mather we discard all $W_i$ for which $W_i \cap Z \neq \emptyset$ or $W_i \cap A = \emptyset$, 
  and reindex the remaining $W_i$'s again by the positive integers. By \eqref{eq:condition} we then have $A \subset U \cup \bigcup_{i\in \NBbb^*} W_i$
  and $W_i \subset  V$ for all $i\in \NBbb^*$. Next, choose $G$-invariant closed subsets $W_i' \subset W_i \cap S$ such that $A \subset U \cup \bigcup_{i\in \NBbb^*} W_i'$.
  Since the $W_x$ are relatively compact, all $W_i'$ are compact. Finally put for all $j \in \NBbb$ 
  \[ U_j = \varphi_0 \big( \pi_{E_0}^{-1} (U) \cap B(\varepsilon_0 ,E_0)\big) \cup W_1 \cup \ldots \cup W_j \ . \]
  Note that the $U_j$ are then $G$-invariant and open and that $U_0 \cap S = U$. 

  We now construct inductively $G$-equivariant diffeotopies $H^0, H^1, H^2, \ldots $ of $M$ together with a sequence $\psi_0, \psi_1,\psi_2,\ldots $
  of $G$-equivariant isomorphims of tubular neighborhoods. We start with defining $H^0_t$ to be the identity map for all $t\in I$ and let $\psi_0$ 
  be the ismorphism from the statement of the theorem. 
  
  For the induction step we assume to be given diffeotopies $H^0,H^1,\ldots , H^{i-1}$ of $M$
  together with $G$-equivariant isomorphisms $\psi_0,\ldots , \psi_{i-1}$ of tubular neighborhoods 
  having the following properties:
  \begin{enumerate}[(a)]
  \item 
    The diffeotopies $H^0,H^1,\ldots , H^{i-1}$  and isomorphisms  $\psi_0,\ldots , \psi_{i-1}$ are 
    $G$-equivariant and compatible with $f$.
  \item 
    The diffeotopies $H^0,H^1,\ldots , H^{i-1}$ leave $S$ pointwise fixed. 
  \item
    For each $j = 0, \ldots , i-1$ the diffeotopy $K^j$ of $M$ defined by 
    $K^j_t := H^j_t \circ H^{j-1}_t \circ \ldots \circ H^0_t$ for $t\in I$ has 
    support in $U_j \cap V$. 
  \item
    One has $\big( K^j_t (x), x\big) \in O$ for all $x \in M$, $t\in I$, and $j = 0, \ldots , i-1$.
  \item
    For each $j = 0, \ldots , i-1$ there exist $G$-invariant relatively compact open neighborhoods $U^*_j$ of
    $A \cup W_1' \cup \ldots \cup W_j'$ in $S$ such that 
    $\overline{U}^*_j \subset U^*_{j-1} \cup W_j$ when $j>0$ 
    and such that $\psi_j$ is an isomorphism of tubular neighborhoods 
    $k^j_* \sfT_0|_{\overline{U}^*_j} \to \sfT_1|_{\overline{U}^*_j}$,
    where $k^j :M \to M$ is the diffeomorphism $K^j_1$.
  \end{enumerate}
  
  By the local $G$-equivariant case from {\it Step 1} there exist a $G$-equivariant diffeotopy $H^i$ on $M$ 
  together with an isomorphism of tubular neighborhoods $\psi_i$ such that the conditions of the induction are 
  satisfied. Let us provide a detailed argument by adapting Mather's argument to the $G$-equivariant case.
  First choose a $G$-equivariant relatively compact open subset  $W^0_i$ of $W_i$ with $W_i' \subset W^0_i$. 
  Then let $U^*_i$ be a $G$-equivariant open neighborhood of $A \cup  W_1' \cup \ldots \cup W_i'$ in $S$
  with closure being compact and in $U^*_{i-1} \cup W^0_i$. By the local $G$-equivariant case there exists
  a diffeotopy $H^j$ of $W_i$ which is $G$-equivariant and compatible with $f$, 
  has support in $W^0_i$ and leaves $S \cap W_i$ invariant. Moreover, since 
  $Z_i := \overline{U}^*_i- U^*_{i-1}$ is a $G$-invariant and compact subset of $W_i$
  and $k^{i-1}_*\sfT_0|_{\overline{U}^*_{i-1}\cap W_i} \sim \sfT_1|_{\overline{U}^*_{i-1}\cap W_i}$, 
  the diffeotopy $H^i$ can be chosen so that there exists a $G$-equivariant isomorphism of tubular neighborhoods 
  \[
   \psi_i : (H^i_1)_*k^{i-1}_*\sfT_0|_{\overline{U}^*_i \cap W_i } \to \sfT_1|_{\overline{U}^*_i \cap W_i} 
  \]
  which fulfills 
  $\psi_i|_{\overline{U}^*_i \cap W_i \cap\overline{U}^*_{i-1}} = \psi_{i-1}|_{\overline{U}^*_i \cap W_i \cap\overline{U}^*_{i-1}}$. 
  Finally, one can even achieve that the $H^i_t$ with $t\in I$ are arbitrarily and uniformly close to the 
  identity. 
  Since the support of the diffeotopy $H^i$ is a compact $G$-invariant subset of $W_i$, one can extend 
  $H^i$ by the identity outside $W_i$ to a $G$-invariant diffeotopy on $M$ which has support in $W_i$ and is compatible with $f$. By putting 
  $\psi_i|_{\overline{U}^*_{i-1}} = \psi_{i-1}|_{\overline{U}^*_{i-1}}$, the isomorphism $\psi_i$ can be extended to 
  $\overline{U}^*_i$ and the thus extended isomorphism has all the desired properties. This completes the induction step. 

  Since the Lie group $G$ is compact, one can shrink the $G$-invariant open neighborhood $O \subset M \times M$ of the diagonal  
  so that the projection $\operatorname{pr}_2 : \overline{O} \to M$ onto the second factor is proper. 
  The sequences $\big( K^i_t  \big)_{i\in \NBbb}$ and  $\big( \psi_i  \big)_{i\in \NBbb}$ then eventually become locally constant.
  Hence the maps 
  \[
    H : \: M \times I \to M, \quad (x,t) \mapsto \lim_{i\to \infty} K^i_t(x)
    \quad \text{and}\quad \psi : \: M \to \Hom (E_0, E_1) , \quad x \mapsto \lim_{i\to \infty} \psi_i(x)
  \]
  By construction, $H$ then is a $G$-equivariant diffeotopy of $M$ compatible with $f$ and $\psi : (H_1)_*\sfT_0 \sim \sfT_1$ 
  a  $G$-equivariant isomorphism of tubular neighborhoods compatible with $f$ as well. Moreover, $H$ and $\psi$ have
  the properties claimed in the theorem. 
\end{proof}

\begin{theorem}[Existence of equivariant tubular neighborhoods]\label{thm:tubular-existence}
   Let $M,N$ be $G$-manifolds, $S \subset M$ a $G$-invariant smooth submanifold, and $f: M\to N$ a
  $G$-equivariant smooth map which is submersive over $S$.  Let $U \subset S$ be  relatively open $G$-invariant subset, 
  and $A \subset U$ relatively closed and $G$-invariant. Assume that $\sfT_0$ is a $G$-equivariant tubular neighborhood
  of $U$ in $M$ compatible with $f|_{T_0}$. Then there exists a $G$-equivariant tubular neighborhood $\sfT$ of $S$
  compatible with $f$ such that $\sfT|_A$ and $\sfT_0|_A$ are $G$-equivariantly isomorphic.  
\end{theorem}

\begin{proof}
  {\it Step 1.}
  The Equivariant Submersion Theorem entails existence of tubular neighborhoods in the local equivariant case.
  Let us explain this. The global case will be considered in the following step.  Assume that $M$ is
  of the form $G \times_H (B \times C \times D)$, $P$ is equivariantly diffeomorphic to $G\times_KB$, and under these identifications
  $S$ has the form $G\times_H (B \times C)$ and $f$ the form $\overline{\id_G \times \Pi}$. 
  Hereby, $H\subset K \subset G$ are closed subgroups, $B \subset \R^p$, $C\subset \R^k$, and 
  $D\subset \R^l$ are open convex neighborhoods of the origin, where $\R^p$ carries an orthogonal $K$-repersentation, 
  and $\R^k$, $\R^l$ orthogonal $H$-representations, and $\Pi$ is projection onto the third factor. 
  Now let $E$ be the bundle $G\times_H (B \times C \times \R^l) \to S \cong G\times_H (B\times C)$,
  $\varepsilon : S \to (0,\infty)$ a constant map such that the ball of radius $\varepsilon$ in $\R^l$ 
  is contained in $D$,  and $\varphi : B(\varepsilon,E) \hookrightarrow M$ the identical embedding.
  Then $\sfT =(E,\varepsilon,\varphi)$ is a $G$-equivariant tubular neighborhood compatible with $f$.

  {\it Step 2.}
  We adapt Mather's argument in the proof of \cite[Prop.~6.2]{MatNTS} to the equivariant case. Without loss of generality we can assume that
  $S$ is closed in $M$.
  Now choose $G$-invariant relatively compact open neighborhoods $W_i$, $i\in \NBbb^*$ together with equivariant charts
  $\Phi_i : W_i \hookrightarrow G \times_{G_i} \R^{p+k+l}$ and $\Psi_i : W_i \hookrightarrow G \times_{G_i} \R^{p+k+l}$  
  fulfilling conditions \eqref{cond:one} to \eqref{cond:three} in  {\it Step 3} of the preceding proof 
  such that the family $(W_i)_{i\in \NBbb^*}$ is a locally finite covering of $S$. Next choose $G$-invariant closed subsets $W_i'\subset S \cap W_i$ 
  such that the family $(W_i')_{i\in \NBbb^*}$ covers $S$ as well. 
  Put $U_0 := T_0 = \varphi_0 \big( B(\varepsilon_0, E_0 \big)$ and define inductively $U_i := W_i \cup U_{i-1}$ for $i\in \NBbb^*$.
  Furthermore put $U_0' := A$ and $U_i' := W_i' \cup U_{i-1}'$ for $i\in \NBbb^*$. Finally let $U_0''$ be a $G$-invariant relatively open neighborhhood 
  of $A$ in $S$ such that $\overline{U}_i'' \subset U$ and then choose  inductively for all $i \in \NBbb^*$ relatively open neighborhoods 
  $U_i''$ of $U_i'$ in $S$ such that $U_i''$ is contained in $W_i\cup U_{i-1}''$ and such that $U_i''$ can be decomposed into
  $G$-invariant relatively open subsets $X_i,Y_i \subset U_i''$ so that $\overline{X}_i \subset W_i \setminus  U_{i-1}'$, 
  $\overline{Y}_i \subset U_{i-1}'' $ and so that $X_i \cap Y_i$ is relatively compact in  $W_i$.
 
  Now we inductively construct $G$-equivariant tubular neighborhoods $\sfT_i$ of $U_i''$ in $M$. The tubular neighborhood $\sfT_0$ is the given one.
  Assume that for some $i\in \NBbb^*$ a $G$-equivariant tubular neighborhood $\sfT_{i-1}$ of $U_{i-1}''$ in $M$ has been constructed and that 
  it is compatible with $f$.  
  By {\it Step 1} there exists a $G$-equivariant tubular neighborhood $\sfT_i'$ of $W_i\cap S$ in $W_i$ which is compatible with $f$.
  So we have two $G$-equivariant tubular neighborhoods over the $G$-invariant subset $ U_i''\cap W_i \cap S$, the corresponding restrictions of 
  $\sfT_{i-1}$ and $\sfT_i'$.
  By Theorem \ref{thm:uniqueness-equivariant-tubular-neighborhoods} there exists a $G$-quivariant diffeomorphism $h$ of $M$
  which is compatible with $f$ and has support within a sufficiently small relatively compact neighborhood of $\overline{X_i \cap Y_i}$
  such that $h_* \sfT_{i-1}|_{X_i \cap Y_i} =  \sfT_i'|_{X_i \cap Y_i}  $. By $h$ having a sufficiently small support we in particular mean that 
  $h$ is the identity in a neighborhood of $U_{i-1}'$. One can now glue together $h_* \sfT_{i-1}$ and $\sfT_i'$ to a $G$-equivariant
  tubular neighborhood $\sfT_i$ over $U_i'' = X_i \cup Y_i$. By construction $\sfT_i$ is compatible with $f$.
  
  Since for all $i\in \NBbb^*$ the tubular neighborhoods $\sfT_{i-1}$ and $\sfT_i$ are isomorphic over a small neighborhood of
  $U_{i-1}'$ in $S$ there exists a $G$-equivariant tubular neighborhood $\sfT$ of $S$ in $M$ such that
  $\sfT|_{U_i'} \sim \sfT_i|_{U_i'}$ for all $i$. This tubular neighborhood is compatible with $f$ since all the $\sfT_i$ are 
  and fulfills the claim. The theorem is proved.      
\end{proof}
\subsection{Existence of equivariant control data}

Before proving the existence of $G$-equivariant control data in Theorem \ref{thm:control-data-existence} below, we first need the following equivariant analog of \cite[Lem. 7.3]{MatNTS}. Given a stratum $S$, a tubular neighbourhood $\mathsf{T} = (E, \varepsilon, \varphi)$ and a smooth function $\varepsilon' : S \rightarrow \mathbb{R}_{>0}$, define $T_{\varepsilon'}^\circ:= \varphi(B_\varepsilon \cap B_{\varepsilon'})$. 

\begin{lemma}\label{lem:control-submersion}
Let $R$ and $S$ be disjoint submanifolds of $M$ which are preserved by $G$, such that the pair $(S,R)$ satisfies condition (B). Let $\mathsf{T}$ be a $G$-equivariant tubular neighbourhood of $R$ in $M$. Then there exists a $G$-invariant smooth function $\varepsilon' : R \rightarrow \mathbb{R}_{>0}$ such that the mapping
\begin{equation*}
(\pi_T, \varrho_T) : S \cap T_{\varepsilon'}^\circ \rightarrow R \times (0, \infty) 
\end{equation*}
is a smooth submersion.
\end{lemma}

\begin{proof}
Since $G$ is compact, then the result follows from the non-equivariant version in \cite[Lem. 7.3]{MatNTS} by averaging over the $G$-orbits in $R$.
\end{proof}

\begin{theorem}\label{thm:control-data-existence}
  Let $G$ be a compact Lie group and $M,N$  smooth $G$-manifolds. Assume that $(X,\calS)$ is a 
  (B) regular stratified subspace of $M$, that $X$ is invariant under the $G$-action 
  and that the induced $G$-action on $X$ is compatible with the stratification $\calS$.
  Assume further that  $f:X \to N$ is a $G$-equivariant smooth stratified submersion. Then there exists a system of 
  $G$-equivariant control data $\calT = (T_S,\pi_S,\varrho_S)_{S\in \calS}$ on $(X,\calS)$ compatible with $f$.     
\end{theorem}

\begin{proof}
The proof is by induction on the dimension of the strata, following the strategy of \cite[Sec. 7]{MatNTS}. Let $\calS_k$ be the subset of $\calS$ consisting of strata of dimension less than or equal to $k$, and let $X_k$ be the union of all strata in $\calS_k$.

Since the strata in $\calS_0$ all have dimension zero, then there exists a system of $G$-equivariant control data $\calT_0 = (T_S, \pi_S, \varrho_S)_{S \in \calS_0}$ on $(X_0, \calS_0)$ which is compatible with $\left. f \right|_{X_0}$.

Now suppose that there exists a system of $G$-equivariant control data $\mathcal{T}_{k-1} = (T_S, \pi_S, \varrho_S)_{S \in S_{k-1}}$ on $(X_{k-1}, \mathcal{S}_{k-1})$ which is compatible with $\left. f \right|_{X_{k-1}}$.

Let $S$ be a stratum of dimension $k$, and for each $\ell = 0, \ldots, k$, define
\begin{equation*}
U_\ell := \bigcup_{Y < S, \dim Y \geq \ell} T_Y , \quad S_\ell := U_\ell \cap S .
\end{equation*}

For each $\ell$, we will construct a tubular neighbourhood $T_\ell$ of $S_\ell$ which satisfies the equivariant control data relations (CC1)--(CC4). Using the approach of \cite[Proof of Prop. 7.1]{MatNTS}, we will do this by descending induction on $\ell$. Note that it is sufficient to construct $T_\ell$ separately for each stratum $Y$ of dimension $\ell$, since if $Y, Y'$ both have dimension $\ell$ then $T_Y \cap T_{Y'} = \emptyset$. 

For the base case $\ell = k$, note that $S_k = \emptyset$, and so there is nothing to prove.

Now suppose that we have constructed $T_{\ell+1}$ such that $\varrho_{\ell+1} : T_{\ell+1} \rightarrow \mathbb{R}$ is $G$-invariant, $\pi_{\ell+1} : T_{\ell+1} \rightarrow S_{\ell+1}$ is $G$-equivariant, and if $Y < S$, $\dim Y \geq \ell+1$, $m \in T_{\ell+1} \cap T_Y$, then
\begin{align}\label{eqn:control-relations-ell+1}
\begin{split}
\varrho_Y \circ \pi_{\ell+1}(x) & = \varrho_Y(x) \\
\pi_Y \circ \pi_{\ell+1}(x) & = \pi_Y(x) .
\end{split}
\end{align}
If necessary, shrink the neighbourhood $T_{\ell+1}$ so that $x \in T_{\ell+1}$ implies that there exists a stratum $Z < S$ with $\dim Z \geq \ell+1$ such that if $x$ is also in $T_Z$ then $\pi_{\ell+1}(x) \in T_Z$.

Given $x \in T_{\ell+1} \cap T_Y$ such that $\pi_{\ell+1}(x) \in T_Y$, then there exists $Z < S$ with $\dim Z \geq \ell+1$, $x \in T_Z$ and $\pi_{\ell+1}(x) \in T_Z$. Therefore $\pi_{\ell+1}(x) \cap T_Y \cap T_Z$ and so $T_Y \cap T_Z$ is non-empty, hence $Y < Z$. Note that the relations (CC1)--(CC4) hold for the pair $(Y, Z)$ by the inductive hypothesis, and also that since $\dim Z \geq \ell+1$ then \eqref{eqn:control-relations-ell+1} holds with $Y$ replaced by $Z$. Therefore we have
\begin{align*}
\varrho_Y \circ \pi_{\ell+1}(x) & = \varrho_Y \circ \pi_Z \circ \pi_{\ell+1}(x) = \varrho_Y \circ \pi_Z(x) = \varrho_Y(x) \\
 \pi_Y \circ \pi_{\ell+1}(x) & = \pi_Y \circ \pi_Z \circ \pi_{\ell+1}(x) = \pi_Y \circ \pi_Z(x) = \pi_Y(x)  .
\end{align*}
Again, since $\dim Y < k$, then we can further suppose from (CC3) that $(\rho_Y, \pi_Y) : T_Y \cap S \rightarrow \mathbb{R} \times Y$ is a submersion, and from (CC1) that $\varrho_Y$ is $G$-invariant and $\pi_Y$ is $G$-equivariant.

Therefore we have constructed a tubular neighbourhood $T_{\ell+1} \cap T_Y \rightarrow S_{\ell+1} \cap T_Y$ and so it only remains to extend it to a neighbourhood $T_{S, Y} \rightarrow S \cap T_Y$ and then to a neighbourhood $T_S \rightarrow S$.

%{\bf Now extend using the tubular neighbourhood theorem as in Mather.}

Now if $S_{\ell+1}^\circ$ is an open subset of $S$ whose closure lies in $S_{\ell+1}$, then Theorem \ref{thm:tubular-existence} shows that there exists a tubular neighbourhood $T_{S, Y}$ of $T_Y \cap S$ such that
\begin{align*}
\varrho_Y \circ \pi_{S, Y}(x) &  = \varrho_Y(x) \\
\pi_Y \circ \pi_{S, Y}(x) & = \pi_Y(x) ,
\end{align*}
the map $\pi_{S, Y}$ is $G$-equivariant and the function $\varrho_{S, Y}$ is $G$-invariant, and such that the restriction of $T_{S, Y}$ to $|T_Y| \cap S_{\ell+1}^\circ$ is isomorphic to the restriction of $T_{\ell+1}$.

Now in the same way as the second step of \cite[Proof of Prop. 7.1]{MatNTS}, we can inductively extend the tubular neighbourhood to a neighbourhood $T_S$ of all of $S$, which is compatible with the submersion $f$, where we use Theorem \ref{thm:tubular-existence} and Lemma \ref{lem:control-submersion} in place of \cite[Prop. 6.2 \& Lem. 7.3]{MatNTS} in order to guarantee that the tubular neighbourhoods are $G$-equivariant. This completes the inductive step, and hence also the proof of the theorem.
\end{proof}

%% file: Equivariant_Wall_result.tex
\section{A stratification compatible with a given set of subvarieties}\label{sec:equivariant-Wall}

In this section we use a construction due to Wall \cite{WalRS} to prove the following theorem.

\begin{theorem}\label{thm:equivariant-Wall}
Let $G$ be a Lie group acting smoothly on a real analytic manifold $M$, and let $( A_r )_{r=1}^n$ be a finite family of analytic subvarieties, each of which is preserved by the action of $G$. Then there is a (B) regular stratification of $M$ in which each $A_r$ is a finite union of $G$-invariant strata.
\end{theorem}

%\begin{remark}
%A result of Kempf \cite[Lem. 1.1]{KempfIIT} shows that the action of any reductive group $G$ on an affine variety $M$ can be linearised, i.e. there is a linear representation space $V$ for $G$ and a $G$-equivariant isomorphism between $M$ and an affine subvariety of $V$. Therefore we do not lose any generality in the above theorem by assuming that the action of $G$ on $M$ is linear.
%\end{remark}

As a preparation for the proof of Theorem \ref{thm:equivariant-Wall}, we make the following observation and prove some preliminary results.

\begin{remark}
 Recall that by the solution of (a variant of) Hilbert's Fifth Problem by  Matumoto–Shiota \cite{MatShiUTOSDTGA} any smooth manifold $M$ with a
 smooth compact Lie group action carries a real-analytic structure so that the Lie group action becomes analytic. 
 Together with the Morrey--Grauert Theorem \cite{MorAEARAM,GraLPIRAM} which tells that every smooth manifold has a unique real-analytic structure this means
 that it is no loss of generality when we assume that the action of the compact Lie group $G$ on the real-analytic manifold $M$ is analytic. 
\end{remark}

\begin{lemma}\label{lem:singular-set-invariant}
Let $M$ be a real analytic manifold, and let $G$ be a connected Lie group acting analytically on $M$. 
If $X$ is a subvariety of $M$ preserved by $G$, then the singular set $X_\textup{sing}$ is also preserved by $G$.
\end{lemma}

\begin{proof}
Given $p \in X$, let $f_1, \ldots, f_n $ be real analytic functions defining $X$ in a neighbourhood of $p$. Any $g \in G$ defines a diffeomorphism 
$\psi_g : M \rightarrow M$. In particular, since $X$ is preserved by the action of $G$ then $ f_1 \circ \psi_g^{-1}, \ldots, f_n \circ \psi_g^{-1} $ define $X$ 
in a neighbourhood of $g \cdot p$. The Jacobian of these equations is $df \circ d \psi_g^{-1}$, which has the same rank as $df$. 
Therefore $g \cdot p$ is a singular point if and only if $p$ is singular, and so $X_\textup{sing}$ is preserved by the action of $G$. 
\end{proof}

\begin{lemma}\label{lem:closure-invariant}
Let $M$ be a metrizable topological space, and let $G$ be a group acting continuously on $M$. If $X \subset M$ is any subset preserved by the action of $G$, 
then $\overline{X}$ and $\overline{X} \setminus X$ are also preserved by the action of $G$.
\end{lemma}

\begin{proof}
Given $p \in \overline{X} \setminus X$, let $( p_n )_{n \in \NBbb} \subset X$ be a sequence in $X$ converging to $p$. Since the action of $G$ is continuous, then for any $g \in G$ the sequence $( g \cdot p_n )_{n\in \NBbb} \subset X$ converges to $g \cdot p$. Since $G$ preserves $X$ and $p \notin X$ then $g \cdot p \notin X$ also. Therefore $g \cdot p \in \overline{X} \setminus X$ for all $g \in G$, and therefore $\overline{X} \setminus X$ is also preserved by the action of $G$, hence so is $\overline{X}$.
\end{proof}

\begin{lemma}\label{lem:regular-invariant}
Let $G$ be a Lie group acting smoothly on a smooth manifold $M$, and let $X$ and $Y$ be two disjoint strata in a stratification of $M$. 
Then $G$ preserves the set of points $x \in X \cap \overline{Y}$ where $(X, Y)$ is (B) regular.
\end{lemma}

\begin{proof}

Given $g \cdot x$, let $\psi_g : M \rightarrow M$ denote the diffeomorphism associated to the action of $g \in G$. Since $G$ acts smoothly on $M$, then for each $g \in G$ a chart $\varphi : U \rightarrow \R^d$ around $x \in M$ determines a chart $\varphi \circ \psi_g^{-1} : g(U) \rightarrow \R^d$ around $g \cdot x \in M$. Since Whitney's condition (B) is independent of the choice of chart (cf. \cite[Lem. 1.4.4]{PflAGSSS}) then $(X, Y)$ is (B) regular at $x \in X$ if and only if $(X, Y)$ is (B) regular at $g \cdot x$.
\end{proof}

We can now use the above results to prove the main theorem of the section.

\begin{proof}[Proof of Theorem \ref{thm:equivariant-Wall}]
We closely follow the proof of the corresponding result of Wall \cite{WalRS} when the group action is trivial, and use the above results to show that the construction extends to the equivariant setting.

Suppose that there exists a filtration $T_i \subset T_{i+1} \subset \cdots \subset T_m = M$ such that
\begin{itemize}

\item each $T_j$ is a closed semianalytic set in $M$,

\item each $T_j$ is preserved by the action of $G$,

\item for each $j = i+1, \ldots, m$, the set $S_j = T_j \setminus T_{j-1}$ is a $j$-dimensional real-analytic manifold
      called the $j$-the stratum, and 
\item each $A_r \cap S_j$ is a union of components of $S_j$.

\end{itemize}
The above conditions are clearly satisfied for $T_m = M$, thus giving us the base case for the induction. Define
\begin{equation*}
B_1 = \begin{cases} (T_i)_\textup{sing} & \text{if $\dim T_i = i$} \\ T_i & \text{if $\dim T_i < i$} \end{cases} 
\end{equation*}
and
\begin{equation*}
B_2 = \bigcup_r \left( A_r \cap \overline{(T_i)_\textup{reg} \setminus A_r} \right) ,
\end{equation*}
By construction $B_1$ and $B_2$ are semianalytic. Lemma \ref{lem:singular-set-invariant} shows that $B_1$ is preserved by $G$, therefore so is $(T_i)_{reg}$. 
Together with Lemma \ref{lem:closure-invariant} this implies that $B_2$ is also preserved by $G$. Hence $T_i \setminus (B_1 \cup B_2)$ is preserved by $G$.
To finish the proof, define the set $B_3$ of points of $T_i$ where some higher-dimensional stratum fails to be (B) regular. 
Lemma \ref{lem:regular-invariant} shows that this is preserved by $G$. By \cite[ p.~337, Proposition]{WalRS}, $B_3$ is semianalytic of dimension less 
than the dimension of $T_i$. 
Then define $T_{i-1} := B_1 \cup B_2 \cup B_3$ and $S_i = T_i \setminus T_{i-1}$. By construction and the above arguments these sets are both 
semianalytic and $G$-invariant. Since $T_i \setminus B_1$ either coincides with $(T_i)_\textup{sing}$ or is empty, and $\dim (B_2 \cup B_3) < i$,
the stratum $S_i$ is a real analytic manifold. 
Moreover, as in Wall \cite{WalRS} one argues that for each $r$ the intersection $A_r \cap S_i$ has no relative frontier by construction and thus 
is a union of components. 
So we can continue inductively to define a (B) regular $G$-invariant stratification by real analytic manifolds such that each $A_r$ is a finite union of strata.
\end{proof}

%% file: Equivariant_deformation_retract.tex
\section{Constructing the equivariant neighbourhood deformation retract}\label{sec:equivariant-def-retract}

Let $M$ be a smooth manifold equipped with the action of a compact Lie group $G$, and let $A \subset X \subset M$ be closed subsets with inclusion map denoted $i : A \hookrightarrow M$. Suppose that $X$ carries a (B) regular $G$-invariant Whitney stratification $\{ S \}_{S \in \mathcal{S}}$, which restricts to a (B) regular $G$-invariant Whitney stratification $\{ S \}_{S \in \mathcal{S}_A}$ of $A$. Theorem \ref{thm:control-data-existence} shows that there exists a system of $G$-equivariant control data on $(X, \mathcal{S})$ and Theorem \ref{thm:equivariant-Wall} shows that these assumptions are satisfied when $A$ and $X$ are $G$-invariant analytic subvarieties of $M$ with $A \subset X$.

In this section we prove Theorem \ref{thm:equivariant-stratified-Verona} which shows that the inclusion $A \hookrightarrow X$ is an equivariant cofibration of stratified spaces. In particular, the result of Corollary \ref{cor:equivariant-stratified-Verona} shows that the homotopy equivalences in the Morse theory of \cite{WilEMTNSMMV} can be chosen to be $G$-equivariant.

Using Theorem \ref{thm:control-data-existence}, construct a system of $G$-equivariant control data $(T_S, \pi_S, \rho_S)_{S \in \mathcal{S}}$ for $X$. Since $A$ is a $G$-invariant stratified subspace of $X$ then $(T_S, \pi_S, \rho_S)_{S \in \mathcal{S}_A}$ is a system of $G$-equivariant control data for $A$. On restricting to a small enough open neighbourhood of $A$, we can assume that
\begin{enumerate}

\item if $S \subset X$ is a stratum of lowest dimension, then $S \subset A$, and 

\item if $S \subset X$ is a stratum of $X$ then $\bar{S} \cap A \neq \emptyset$.

\end{enumerate}

First we set up some notation and prove some preliminary results. On each tubular neighbourhood $T_S$, fix a radial vector field $\frac{\partial}{\partial \rho_S}$ as in \cite[Cor. 3.7.4]{PflAGSSS}. Since $\rho_S$ is $G$-invariant then $\frac{\partial}{\partial \rho_S}$ is $G$-equivariant and so is its integral flow. Using the integral flow of radial vector fields, for each stratum $S$ and each $x \in S$, there exists a neighbourhood $U_x \subset T_S$ and a real number $r > 0$ together with an isomorphism of stratified spaces 
\begin{equation}\label{eqn:mapping-cylinder}
U_x \cong (\rho_S^{-1}(r) \cap U_x) \times [0,r] / \sim ,
\end{equation}
where $(y_1, 0) \sim (y_2, 0)$ if and only if $\pi_S(y_1) = \pi_S(y_2)$. Equivalently, $U_x$ is homeomorphic to the mapping cylinder of $\left. \pi_S \right|_{U_x \cap \rho_S^{-1}(r)}$ and this homeomorphism is determined by the flow of the radial vector field $\frac{\partial}{\partial \rho_S}$.

Given any stratum $S_0 \in \mathcal{S}_A$ and a sequence of strata $S_0 < S_1 < \cdots < S_k \subset A$ of increasing height, define
\begin{equation}\label{eqn:intersection-neighbourhood}
T_{S_0, \ldots, S_k} := T_{S_0} \cap T_{S_1} \cap \cdots \cap T_{S_k} \setminus \left( S_0 \cup S_1 \cup \cdots \cup S_k \right) .
\end{equation}
Given any $x \in T_{S_0, \ldots, S_k}$, there exists a neighbourhood $U$ of $x$ such that $U$ is contained in a trivialisation for each of $\pi_{S_0}$, $\pi_{S_1}$, \ldots, $\pi_{S_k}$. Therefore there exist $r_0, \ldots, r_k$ and $\varepsilon > 0$ such that
\begin{equation*}
V := \bigcap_{\ell=0}^k \rho_{S_\ell}^{-1}((r_\ell - \varepsilon, r_\ell + \varepsilon)) \subset U .
\end{equation*}
Let $Y := \bigcap_{\ell=0}^k \rho_{S_\ell}^{-1}(r_\ell)$. Again using the integral flow of radial vector fields, we have
\begin{equation}\label{eqn:product-neighbourhood}
V \cong Y \times \prod_{\ell=0}^k (r_\ell - \varepsilon, r_\ell + \varepsilon) ,
\end{equation}
and for any $x = (y, t_0, \ldots, t_k) \in V$ we have $\rho_{S_\ell}(x) = t_\ell$ for each $\ell=0, \ldots, k$.

For each stratum $S_\ell$, let $\varphi_{S_\ell}$ be the integral flow of the radial vector field on the tubular neighbourhoood $T_{S_\ell}$. Recall that these flows have the following properties
\begin{itemize}

\item $\varphi_{S_\ell}$ preserves the tubular distance functions $\rho_{S_j}$ for each $j \neq \ell$, 

\item $\varphi_{S_\ell}$ is $G$-equivariant,

\item the flow on the cylinder $\rho_{S_\ell}^{-1}(r_\ell) \times (r_\ell - \varepsilon, r_\ell + \varepsilon)$ is given by $\varphi_{S_\ell} ((y_\ell, t_\ell), t) = (y_\ell, t_\ell+t)$, and

\item the flow preserves strata.

\end{itemize}

Therefore, on the neighbourhood $V = \bigcap_{\ell=0}^k \rho_{S_\ell}^{-1}((r_\ell - \varepsilon, r_\ell + \varepsilon)) \cong Y \times \prod_{\ell=0}^k (r_\ell - \varepsilon, r_\ell + \varepsilon)$ the flow is given by
\begin{equation*}
\varphi_{S_\ell}((y,t_0, \ldots, t_k), t) = (y, t_0, \ldots, t_\ell + t, \ldots, t_k) .
\end{equation*}
In particular, $\varphi_{S_\ell}$ is the flow of the vector field $\frac{\partial}{\partial t_\ell}$ on $Y \times \prod_{\ell=0}^k (r_\ell - \varepsilon, r_\ell + \varepsilon)$, the vector fields $\{ \frac{\partial}{\partial t_\ell} \}_{\ell=0, \ldots, k}$ are linearly independent, and the flows $\varphi_{S_\ell}$ for $\ell = 0, \ldots, k$ all commute and preserve strata. Moreover, even though the above calculations have been done with respect to the local neighbourhood $V$, these vector fields and flows are well-defined and $G$-equivariant on the entire neighbourhood $T_{S_0, \ldots, S_k}$, since the radial vector fields are well-defined and $G$-equivariant on $T_{S_0, \ldots, S_k}$. 

Given functions $a_\ell : V \rightarrow \R$ for each $\ell = 0, \ldots, k$, the vector field 
\begin{equation*}
\chi(y, t_0, \ldots, t_k) = \sum_{\ell = 0}^k a_\ell(y, t_0, \ldots, t_k) \frac{\partial}{\partial t_\ell}  
\end{equation*}
is also tangent to strata, and so the flow preserves strata. Moreover, if the functions $a_\ell$ are independent of $y$, then this vector field is $G$-equivariant and hence the flow is $G$-equivariant, since the $G$-action preserves the radial distance functions $\rho_{S_\ell}$.

The next lemma is used in the proof of Theorem \ref{thm:equivariant-stratified-Verona}. 

\begin{lemma}
Let $Q = [-1,1]$, $B = [0,1] \times [0,1]$ and $C = (\{0\} \times [0,1]) \cup ([0,1] \times \{0\}) \subset B$. Then there exists a proper continuous mapping $H : Q \times [0,1] \rightarrow B$ such that
\begin{equation*}
H(Q \times (0,1)) \subset B \, \setminus \, C, \quad H(Q \times \{0\}) = C 
\end{equation*}
and $\left. H \right|_{Q \times (0,1)}$ is a diffeomorphism onto its image.
\end{lemma}

\begin{proof}
Choose a smooth monotone function $\phi : [0, \frac{\pi}{2}] \rightarrow [0, \pi]$ such that $\phi(\theta) = \theta$ if $0 \leq \theta \leq \frac{\pi}{3}$ and $\phi(\theta) = \theta + \frac{\pi}{2}$ if $\frac{2\pi}{3} \leq \theta \leq \frac{\pi}{2}$. For notation, let $P : \{ (x, y) \in \R^2 \mid y \geq 0, (x,y) \neq (0,0) \} \rightarrow \R_{> 0} \times [0,\pi]$ be the polar coordinate homeomorphism. Then the map $h : B \rightarrow [-1,1] \times \R_{\geq 0}$ given by
\begin{equation*}
h \circ P^{-1}(r, \theta) := P^{-1}(r, \phi(\theta)) , \quad h(0,0) = (0,0)
\end{equation*}
is a homeomorphism onto its image, which restricts to a diffeomorphism of $B \, \setminus \, C$ onto $h(B \, \setminus \, C)$. Moreover (in Cartesian coordinates), the image of $h$ contains $[-1,1] \times [0, \frac{1}{2}]$. Now define $H : Q \times [0,1] \rightarrow B$ by $H(q, t) = h^{-1}(q, \frac{t}{2})$.
\end{proof}

Let $W \subset B$ be the image of $\left. H \right|_{Q \times [0,1)}$. The previous lemma shows that $H$ restricts to a diffeomorphism $Q \times (0,1) \cong W \setminus C$. Using the homeomorphism $H$, for any $w \in W$ we can write $w = H(q(w), s(w))$, where $(q(w), s(w)) \in Q \times (0,1)$. Define a flow $\varphi: W \times [0,\infty) \rightarrow W$ by
\begin{equation*}
\varphi(w,t) = \begin{cases} H(q(w), e^{-t} s(w)) & w \notin C \\ w & w \in C \end{cases}
\end{equation*}
Taking the vector field associated to this flow gives us the following lemma.

\begin{lemma}\label{lem:smoothing-corner}
There exist non-negative smooth functions $a, b : W \rightarrow \R_{\geq 0}$ such that the vector field 
\begin{equation*}
X(x,y) = - a(x,y) \frac{\partial}{\partial x} - b(x,y) \frac{\partial}{\partial y}
\end{equation*}
defined on $W$ satisfies the boundary conditions $X(x,0) = 0 = X(0,y)$,
\begin{equation}\label{eqn:X-boundary-condition}
X(x,1) = - x \frac{\partial}{\partial x} \quad \text{for all $x \in \left[0, \tfrac{1}{2}\right]$}, \quad X(1,y) = - y \frac{\partial}{\partial y} \quad \text{for all $y \in \left[0, \tfrac{1}{2}\right]$},
\end{equation}
and the flow of $X$ defines a smooth map 
\begin{equation*}
\varphi : W \times [0, \infty) \rightarrow W
\end{equation*}
such that $\lim_{t \rightarrow \infty} \varphi((x,y), t) \in C$ for all $(x, y) \in W$.
\end{lemma}

Now define the sets
\begin{align*}
W_{1/2} & := H(Q \times \{ \tfrac{1}{2} \}) \cong Q \\
W_{\leq 1/2} & := H(Q \times [0, \tfrac{1}{2}]) \\
W_{<1/2} & := W_{\leq 1/2} \setminus W_{1/2} .
\end{align*}
Note that the flow $\varphi$ of the vector field $X$ from the Lemma \ref{lem:smoothing-corner} defines a deformation retract of $W_{\leq 1/2}$ onto $C$. Moreover, given such a vector field, for any $w \in W \setminus C$ there exists a unique $t = t(w) \in \mathbb{R}$ such that $\varphi(w, -t(w)) \in W_{1/2}$. 

\begin{definition}\label{def:smoothing-radial-distance}
Given $\varepsilon_1, \varepsilon_2  > 0$, identify $W_{1/2} \cong Q \cong [-1,1]$ and choose a smooth monotone function $f : [-1,1] \rightarrow \R$ such that $f(-1) = \varepsilon_1$ and $f(1) = \varepsilon_2$. The \emph{modified radial distance} $\tilde{\rho} : W \rightarrow [0,1]$ is given by 
\begin{equation*}
\tilde{\rho}(x,y) = \begin{cases} e^{-t(w)} f(\varphi(w, -t(w))) & \text{if $w \in W \setminus C$} \\ 0 & \text{if $w \in C$.} \end{cases}
\end{equation*} 
\end{definition}

Now let $h$ be the maximal height of a stratum in $A$. For each $\ell = 0, \ldots, h$, let $\hat{S}_\ell \subset A$ denote the union of all the strata $S \in \mathcal{S}_A$ such that $\height(S) \leq \ell$. Consider a pair $(U, \varphi_\ell)$ consisting of an open set $U \subset X$ containing $\hat{S}_\ell$ and a flow $\varphi_\ell$ defined on $U$. We say that $(U, \varphi_\ell)$ \emph{has property $(R_\ell)$} if all of the following are satisfied.

\begin{enumerate}

\item \label{item:continuous} $\varphi_\ell$ is continuous.

\item \label{item:limit} $\lim_{t \rightarrow \infty} \varphi_\ell(x, t) \in \hat{S}_\ell$.

\item \label{item:fixed-set} $\varphi_\ell(x, t) = x$ for all $x \in \hat{S}_\ell$.

\item \label{item:preserve-strata} For any stratum $S \in \mathcal{S}$, if $x \in S$ then $\varphi_\ell(x, t) \in S$ for all $t \in [0, \infty)$.

\item \label{item:equivariant} $\varphi_\ell$ is $G$-equivariant.

\item \label{item:mapping-cylinder} For each $r \in [0, 1)$, define $U_r := \varphi_\ell(U, -\log(1-r))$, and define $U_1 := A$. Then $U_r$ satisfies the following conditions

\begin{enumerate}

\item $U_r$ is open in $X$ for all $r \in [0,1)$,

\item $U_r = \bigcup_{s > r} U_s$ and $\overline{U_r} = \bigcap_{s < r} U_s$ for all $r \in (0,1)$.

\end{enumerate}

\end{enumerate}

Note that the condition $\overline{U_r} = \bigcap_{s < r} U_s$ for all $r \in (0,1)$ implies that $\overline{U_r} \subset U_s$ for all $r > s$.

The following theorem is the main result of this section.

\begin{theorem}\label{thm:equivariant-stratified-Verona}
Let $M$ be a smooth manifold equipped with the action of a compact Lie group $G$, and let $A \subset X \subset M$ be closed $G$-invariant subsets with inclusion map denoted $i : A \hookrightarrow M$. Suppose that $X$ carries a $G$-invariant  (B) regular Whitney stratification $\{ S \}_{S \in \mathcal{S}}$ and that there exists a subset $\mathcal{S}_A \subset \mathcal{S}$ such that $A = \bigcup_{S \in \mathcal{S}_A} S$, therefore $\{ S \}_{S \in \mathcal{S}_A}$ is a $G$-invariant (B) regular Whitney stratification of $A$. 

Then there exists a $G$-stratified space $\tilde{A}$, a proper continuous map $\eta : \tilde{A} \rightarrow A$, an open neighbourhood $U$ of $A$ in $X$, and a $G$-equivariant homeomorphism of $U$ onto the mapping cylinder $\psi : U \rightarrow Z_\eta = (\tilde{A} \times [0,1]) / \sim$ such that $\left. \psi \right|_A$ is the identity and $\left. \psi \right|_{A \times (0,1]}$ is a homeomorphism of stratified spaces.

\end{theorem}

\begin{proof}
We first construct a pair $(U, \varphi_h)$ which has property $(R_h)$ by constructing a $G$-equivariant radial vector field, and then define the space $\tilde{A}$ at the end of the proof. 

Consider the neighbourhood $U^{(0)} = \bigcup_{\height(S) = 0} T_S$ of $\hat{S}_0$, and define the vector field $X_0 = -\rho_S \frac{\partial}{\partial \rho_S}$. Note that the vector field is well-defined as the tubular neighbourhoods do not overlap since the strata $S$ all have the same height. Since the radial distance functions $\rho_{S}$ are $G$-invariant and the radial vector field $\frac{\partial}{\partial \rho_S}$ is $G$-equivariant, the  vector field $X_0$ is also $G$-equivariant and so is its flow. It is easy to check the first four conditions of property $(R_\ell)$. Since the flow is continuous and the tubular distance function $\rho_{S}$ is strictly decreasing,  the remaining condition of property $(R_\ell)$ is also satisfied. Note also that $X_0$ commutes with $\frac{\partial}{\partial \rho_{S'}}$ for each stratum $S'$ such that $\height(S') > 0$.

Now suppose that we have a vector field $X_{\ell-1}$ defined on a $G$-invariant neighbourhood $U^{(\ell-1)}$ of $\hat{S}_{\ell-1}$ with $G$-invariant tubular distance function $\tilde{\rho}_{\ell-1}$ and $G$-invariant tubular size function $\tilde{\varepsilon}_{\ell-1}$ such that $U^{(\ell-1)} = \{ \tilde{\rho}_{\ell-1}(x) < \tilde{\varepsilon}_{\ell-1}(x) \}$ and such that the flow $\varphi_{\ell-1}$ of $X_{\ell-1}$ satisfies property $(R_{\ell-1})$. Suppose also that $X_{\ell-1}$ commutes with $\frac{\partial}{\partial \rho_{S'}}$ for each stratum $S'$ such that $\height(S') \geq \ell$. In analogy with the non-equivariant case studied by Verona \cite{VerHPAP} (see also \cite[Sec. 3.9]{PflAGSSS}), we define a $G$-invariant neighbourhood $U^{(\ell)}$ of $\hat{S}_\ell$ and a vector field $X_{\ell}$ satisfying property $(R_\ell)$ by ``smoothing the corner'' using Lemma \ref{lem:smoothing-corner} as follows. First we define $X_\ell$ on $U^{(\ell-1)} \cup \bigcup_{\height(S) = \ell} T_S$ by 

\begin{enumerate}

\item On the subset $U^{(\ell-1)} \, \setminus \, \left( \bigcup_{\height(S) = \ell} T_S \right)$, define $X_\ell = X_{\ell-1}$.

\item For each stratum $S$ with $\height(S) = \ell$, on the subset $\left(T_S \cap \{ \rho_S(x) < \varepsilon_S(x) \} \right) \, \setminus \, U^{(\ell-1)}$ define $X_\ell = -\rho_S(x) \frac{\partial}{\partial \rho_S}$.

\item For each stratum $S$ with $\height(S) = \ell$, on the subset $\left(T_S \cap \{ \rho_S(x) < \varepsilon_S(x) \} \right) \cap U^{(\ell-1)}$ define 
\begin{equation*}
X_\ell(x) = -a\left(\tfrac{\rho_S(x)}{\varepsilon_S(x)}, \tfrac{\tilde{\rho}_{\ell-1}(x)}{\tilde{\varepsilon}_{\ell-1}(x)} \right) \frac{\partial}{\partial \rho_S} + b \left(\tfrac{\rho_S(x)}{\varepsilon_S(x)}, \tfrac{\tilde{\rho}_{\ell-1}(x)}{\tilde{\varepsilon}_{\ell-1}(x)} \right) X_{\ell-1}(x)
\end{equation*}
where the functions $a$ and $b$ are given by Lemma \ref{lem:smoothing-corner}.

\end{enumerate}

Now restrict to the subset
\begin{equation*}
U' := \left( U^{(\ell-1)} \cap \{ \tilde{\rho}_{\ell-1}(x) < \tilde{\varepsilon}_{\ell-1}(x) \}  \right) \cup \bigcup_{\height(S) = \ell} \left(T_S \cap \{ \rho_S(x) < \varepsilon_S(x) \} \right).
\end{equation*}
The result of Lemma \ref{lem:smoothing-corner} shows that the vector field $X_\ell$ is smooth on $U'$. By setting $\varepsilon_1(x) = \varepsilon_S(x)$ and $\varepsilon_2(x) = \tilde{\varepsilon}_{\ell-1}(x)$ (both of which are $G$-invariant), we can glue the modified radial distance $\tilde{\rho}_\ell$ of Definition \ref{def:smoothing-radial-distance} with the radial distance $\tilde{\rho}_{\ell-1}$ on $U^{(\ell-1)} \, \setminus \, \left( \bigcup_{\height(S) = \ell} T_S \right)$ and the radial distance $\rho_S$ on $T_S \, \setminus \, U^{(\ell-1)}$ for each stratum $S$ of height $\ell$. Since $\rho_S$ and $\tilde{\rho}_{\ell-1}$ are both $G$-invariant then this gives us a smooth $G$-invariant radial distance function $\tilde{\rho}_\ell : U' \rightarrow \R_{\geq 0}$, together with a $G$-invariant size function $\tilde{\varepsilon}_\ell : \hat{S}_\ell \rightarrow \R_{> 0}$ such that $\{ \tilde{\rho}_\ell(x) < \tilde{\varepsilon}_\ell(x) \} \subset U'$. Moreover, for each stratum of height $\ell$, on the subset $T_S \, \setminus \, U^{(\ell-1)}$ we have $\tilde{\rho}_\ell = \rho_S$ and $\tilde{\varepsilon}_\ell = \varepsilon_S$, and on the subset $U^{(\ell-1)} \, \setminus \, \bigcup_{\height(S) = \ell} T_S$ we have $\tilde{\rho}_\ell = \tilde{\rho}_{\ell-1}$ and $\tilde{\varepsilon}_\ell = \tilde{\varepsilon}_{\ell-1}$. 

Now define
\begin{equation*}
U^{(\ell)} := \{ x \in U' \, : \, \tilde{\rho}_\ell(x) < \tilde{\varepsilon}_\ell(x) \} .
\end{equation*}
It only remains to verify that the conditions of property $(R_\ell)$ are satisfied. Since the construction of $X_\ell$ only depends on the $G$-invariant functions $\tilde{\rho}_{\ell-1}$, $\tilde{\varepsilon}_{\ell-1}$, $\rho_S$ and $\varepsilon_S$, as well as the $G$-equivariant vector fields $X_{\ell-1}$ and $\frac{\partial}{\partial \rho_S}$ then $X_\ell$ is $G$-equivariant. Since the vector fields $X_{\ell-1}$ and $\frac{\partial}{\partial \rho_{S}}$ commute and their flows preserve strata, the flow of $X_\ell$ also preserves strata. 

Moreover, since $X_\ell$ is constructed from the vector fields $X_{\ell-1}$ and $\frac{\partial}{\partial \rho_S}$ where $\height(S) = \ell$, which commute with $\frac{\partial}{\partial \rho_{S'}}$ for any stratum $S'$ such that $\height(S') > \ell$, and the functions $a\left(\tfrac{\rho_S(x)}{\varepsilon_S(x)}, \tfrac{\tilde{\rho}_{\ell-1}(x)}{\tilde{\varepsilon}_{\ell-1}(x)} \right)$ and $b \left(\tfrac{\rho_S(x)}{\varepsilon_S(x)}, \tfrac{\tilde{\rho}_{\ell-1}(x)}{\tilde{\varepsilon}_{\ell-1}(x)} \right)$ where $\height(S) = \ell$, which are invariant under the flow of $\frac{\partial}{\partial \rho_{S'}}$ for any stratum $S'$ with $\height(S') > \ell$, then $X_\ell$ commutes with $\frac{\partial}{\partial \rho_{S'}}$ for any stratum $S'$ with $\height(S') > \ell$.  

The vector field from Lemma \ref{lem:smoothing-corner} satisfies the remaining conditions \eqref{item:continuous}--\eqref{item:fixed-set} and \eqref{item:mapping-cylinder} of property $(R_\ell)$, hence $X_\ell$ also satisfies these 
conditions. Therefore we can inductively construct a vector field $X_h$ on 
$U^{(h)}$ whose flow $\varphi_h$ has property $(R_h)$.

In the process of the proof, we constructed a $G$-invariant radial distance function $\tilde{\rho}_h : U^{(h)} \rightarrow \mathbb{R}_{\geq 0}$ and a $G$-invariant size function $\tilde{\varepsilon}_h$. Define a rescaled distance function $\rho : U^{(h)} \rightarrow \mathbb{R}_{\geq 0}$ by
\begin{equation*}
\rho(x) = \begin{cases} 0 & x \in A \\ \frac{\tilde{\rho}_h(x)}{\tilde{\varepsilon}_h(x)} & x \notin A \end{cases}
\end{equation*}
By $G$-invariance of $\tilde{\rho}_h$ and $\tilde{\varepsilon}_h$ the
rescaled distance function $\rho$ is  $G$-invariant as well. 
The space $\tilde{A}$ is then defined to be
\begin{equation*}
\tilde{A} := \rho^{-1}(1/2) 
\end{equation*}
and the map $\eta : \tilde{A} \rightarrow A$ is given by taking the limit of the flow $\varphi_h$.
\end{proof}

This immediately gives us the following result, which shows that the main theorem of Morse theory from \cite[Thm. 1.1]{WilEMTNSMMV} can be made to work in the equivariant setting. 

\begin{corollary}\label{cor:equivariant-stratified-Verona}
Let $M$ be a smooth manifold equipped with the action of a compact Lie group $G$, and let $A \subset X \subset M$ be closed subsets with inclusion map denoted $i : A \hookrightarrow M$. Suppose that $X$ carries a $G$-invariant  (B) regular Whitney stratification $\{ S \}_{S \in \mathcal{S}}$, which restricts to a $G$-invariant (B) regular Whitney stratification $\{ S \}_{S \in \mathcal{S}_A}$ of $A$. 

Then there exists a neighbourhood $U$ of $A$ in $X$ and a $G$-equivariant flow $\varphi : U \times [0,1] \rightarrow X$ defining a deformation retract of $U$ onto $A$ such that $U_t := \varphi(U, t)$ satisfies the following conditions

\begin{enumerate}[{\rm (a)}]

\item $U_s$ is open in $X$ for all $s \in [0,1)$,

\item $U_s = \bigcup_{t > s} U_t$ and $\overline{U_s} = \bigcap_{t < s} U_t$ for all $s \in (0,1)$.

\end{enumerate}

\end{corollary}

 Finally in this section we will sketch how to derive an equivariant version of Thom's First Isotopy Lemma 
 \cite[Prop.~11.1]{MatNTS} from  
 the existence of equivariant control data. 
 To this end assume that $X, M$ are as in the corollary, that $P$ is a smooth $G$-manifold, and that $f: M \to P$ is smooth with restriction to $X$ 
 being an equivariant proper controlled submersion. 
 According to the Equivariant Submersion Theorem
 \ref{prop:LocalFormEquivariantSubmersion} the manifold $M$  looks
 locally around a point $x\in X\subset M$ like 
 $G\times_H (B\times C)$, $P$ around $f(x)$ like $G\times_K B$ 
 and $f$ is identified in this representation locally around $X$ with the map 
 $\overline{id_G \times \pi}$. Here we have used the notation from 
 \ref{prop:LocalFormEquivariantSubmersion}. In particular 
 $H$ coincides with the isotropy group $G_x$, $K$ with $G_{f(x)}$ and 
 $\pi : B \times C \to B$ is projection onto the first coordinate.
 One now verifies that the intersection of $X$ with $C$ is a Whitney (B)
 regular stratified space $F$. Hence locally around  $x$ the space
 $X$ is of the form $G \times_H (B \times F)$. One thus obtains the following. 
 
 \begin{theorem}[Equivariant version of Thom's First Isotopy Lemma]
  Let $M,N$ be smooth $G$-manifolds, $X\subset M$ a closed $G$-invariant subset admitting a $G$-invariant Whitney 
  stratification, and $f:M \to P$ a $G$-equivariant smooth map whose restriction to $X$ is a proper stratified submersion. 
  Then the restriction  $f_{|X} : X \to P$ is equivariantly locally trivial
  which means, using notation from above, that  locally around  $x$, the 
  subspace $X$ is of the form $G \times_H (B \times F)$ and the 
   $f_{|X}$ coincides in this local representation with the 
   ``projection'' $G \times_H (B \times F) \to  G_K\times B$.
 \end{theorem}

\begin{remark}
 An equivariant version of Thom's First Isotopy Lemma has been used by Bierstone \cite{BieGPEM} to show an openness theorem
 for equivariant transversality. The verification of the equivariant version of Thom's First Isotopy Lemma 
 has been left to the reader, see \cite[Sec.~9]{BieGPEM}. In this work we presented the missing details  
 which make Mather's machinery work also in the equivariant setting. 
 The main non-trivial steps  hereby have been the proof of the Equivariant Submersion Theorem 
 Prop.~\ref{prop:LocalFormEquivariantSubmersion} and, as a consequence, the existence of equivariant control data 
 in Theorem \ref{thm:control-data-existence}. 
\end{remark}
\begin{remark}
An application of our main result is to Morse theory on singular spaces
carrying a  compact Lie group action (cf.~\cite{WilEMTNSMMV}). Given a real analytic
manifold $M$ with the action of a Lie group $G$, a $G$-invariant
closed analytic variety $Z \subset M$ and an invariant Morse function $f :M \rightarrow \mathbb{R}$
satisfying some additional conditions, \cite[Thm.~1.1]{WilEMTNSMMV} shows that
the main theorem of Morse theory holds in an equivariant sense within this setting. 
In particular this means that for elements $a< b$ in the image of the restriction $f|_Z$ such that
there is one critical value $c$ in between $a$ and $b$, the set
$Z_b = \{x \in Z \mid f(x) \leq b \}$ is homotopy equivalent to the union of $Z_a$ and the unstable set for the critical value $c$.
Moreover this homotopy equivalence can be chosen to be equivariant. The equivariant version of the
main theorem of Morse theory can be applied to the norm-square of a moment map on a (possibly singular)
affine variety (cf.~\cite[Thm. 1.3]{WilEMTNSMMV}). An important example of such a variety is the space of representations of a quiver satisfying a finite set of relations, for which the topological invariants have important applications in representation theory (see for example \cite{Nakajima04}). 
\end{remark}

%% file: Appendix_stratifications.tex
\section{Stratified spaces in the sense of Mather}
\label{sec:stratified-spaces}

In this paper we use stratified space in the sense of Mather \cite{MatSM}.
Let us briefly recall the definition; for further details see  \cite[Sec.~1.2]{PflAGSSS}. 

By a \emph{prestratification} or \emph{decomposition} of a separable locally compact (Hausdorff) space $X$ 
one understands a partition $\calZ$ of $X$ into 
locally closed subspaces $S \subset X$ each carrying the structure of a 
smooth manifold such that the decomposition is locally 
finite and fulfills the \emph{condition of frontier}. 
The latter means that for each pair $R,S \in \calZ$ with the closure 
of $S$ meeting $R$ the relation $R \subset \overline{S}$ holds true. 
The elements of $\calZ$ are called the \emph{pieces} or \emph{strata} of the decomposition.  
If $R,S$ are two strata of $X$ one calls $R$ \emph{incident} to $S$ if $R \subset \overline{S}$ 
and denotes this by $R \leq S$ respectively by $R < S$ if in addition $R$ is not equal to $S$.  

A \emph{stratification} of a locally compact $X$ now is a map $\calS$ which 
assigns to every point $x$ of $X$ a set germ $\calS_x$ at $x$ such that there 
exists for each $x \in X$ an open neighborhood $U$ of $x$ and a decomposition 
$\calZ$ of $U$ with the property that for every point $y$ in $U$ the set 
germ $\calS_y$ coincides with the set germ $[R]_y$ at $y$ of the piece 
$R \in \calZ$ containing $y$. One calls such a decomposition $\calZ$ 
a decomposition \emph{inducing} the stratification $\calS$ over $U$
or a local $\calS$-decomposition around $x$.  

By a \emph{stratified space} we understand a pair $(X,\calS)$ consisting of a separable 
locally compact space $X$ called the \emph{total space} together with a 
stratification $\calS$ on it. In the following $(X,\calS)$ will always denote a stratified space. 

Given an element $x$ of a stratified space   $(X,\calS)$ one defines its \emph{depth} $\depth (x)$
as the maximal number $d$ such that there exist pieces $S_0 , S_1 , \ldots , S_d$ 
of a local $\calS$-decomposition $\calZ$ around $x$ which fulfill
\[
      x \in S_0 < \ldots < S_d \ . 
\]
The depth of $x$ is actually not dependent on a local $\calS$-decomposition $\calZ$ around $x$, see \cite[Lem.~2.1]{MatSM} or 
\cite[Lem.~1.2.5]{PflAGSSS}.
The depth function is locally constant on each stratum of a local decomposition. It allows 
to define a global decomposition of $X$ inducing the stratification $\calS$.
Namely  for each pair of natural numbers $d,m$ let $S_{d,m}$ be the set of 
points $x \in X$ of depth $d$ and for which the dimension of the set germ $\calS_x$ equals 
$m$. Then $S_{d,m}$ is a smooth manifold and the set $\{ S_{d,m} \mid d,m \in \NBbb \}$
is a global decomposition of $X$ inducing $\calS$. It is the coarsest decomposition
with that property, see \cite[Prop.~1.2.7]{PflAGSSS}. We denote this decomposition by the 
symbol $\calS$ also and call its pieces the strata of $(X,\calS)$. We often write $S \in \calS$ to denote 
that $S$ is a stratum of $(X,\calS)$.
The supremum of all depths $\depth (x)$,
where $x$ runs through the points of $X$,  will be called the \emph{depth} of the stratified space $(X,\calS)$. It can be infinite.  Note that the depth is constant on each stratum so it is clear  
what is meant by the \emph{depth} of a stratum.  It is denoted $\depth (S)$.  

Closely related to the depth is the \emph{height} $\height (R)$ of a stratum $R$. It is defined 
as the maximal natural number $h$ such that there exists strata $R_0 < \ldots < R_h$ with 
$R =  R_h$. 

If $(X,\calS)$ and $(Y,\calR)$ denote stratified spaces, a continuous map $f : X \to Y$ 
is called \emph{stratified}, if $f(\calS_x) \subset \calR_{f(x)}$  for all $x \in X$ 
and if the restriction of $f$ to each connected component $S$ of a stratum of $(X,\calS)$
is a smooth map from $S$ to the stratum $R_S$ of $(Y,\calR)$  containing $f(S)$. 
If in addition all the restrictions $f_{|S}: S \to R_S$ are immersions (resp.~submersions),
one calls $f$ a \emph{stratified immersion} (resp.~\emph{stratified submersion}). 

A subspace $A$ of a stratified space $(X,\calS)$ is called a \emph{stratified subspace}
if the map $\calS_A$ which associates to each point $x\in A$ the set germ $A \cap \calS$ 
is a stratification of $A$. In this case $(A,\calS_A)$ becomes a stratified space and
the canonical injection $i : A \hookrightarrow X$ is a stratified immersion. 
If in addition $i$ is a stratified submersion we call $(A,\calS_A)$ a \emph{submersed stratified subspace}.
A subspace $A \subset X$ is a closed submersed stratified subspace  of $(X,\calS)$ if and only if it is a 
union of connected components of strata of $X$.  

Whitney's regularity conditions (A) and (B) play a crucial role in stratification theory in particular in Mather's 
proof of Thom's isotopy lemmata \cite{MatNTS}.
They describe properties how a stratum  of a stratified space  embedded in a smooth manifold $M$ can approach an incident stratum
near its frontier.  Let us recall the Whitney conditions following  \cite[1.4.3]{PflAGSSS}.
A pair $(R,S)$ of smooth submanifolds of $M$ is said to fulfill \emph{Whitney's condition} (A)
at $x \in R$ or that $(R,S)$ is (A) \emph{regular} at $x$ if the following holds. 
\begin{enumerate}[(A)]
\item[(A)]
Let $(y_k)_{k\in \NBbb}$ be a sequence of points of $S$ converging to $x$ 
such that the sequence $T_{y_k}S$, $k \in \NBbb$, of tangent spaces converges in the Gra{\ss}mannian bundle of 
$\dim S$-dimensional subspaces of $TM$ to some $\tau \subset T_xM$. Then $T_xR \subset \tau$.  
\end{enumerate}
The pair $(R,S)$ is said to fulfill \emph{Whitney's condition} (B)
at $x \in R$ or that $(R,S)$ is (B) \emph{regular} at $x$ if for some chart $\chi: U \to \R^d$ of $M$ around $x$
the following is satisfied. 
\begin{enumerate}[(B)]
\item[(B)]
Let $(y_k)_{k\in \NBbb}$ be a sequence in $S$ and $(x_k)_{k\in \NBbb}$ a sequence 
in $R$ such that both converge to $x$ and such that $x_k \neq y_k$ for all $k \in \NBbb$. Assume that the sequence 
of lines $\overline{\chi (x_k)\chi(y_k)}$, where $k$ is large enough so that $x_k, y_k \in U$, converges  in projective space $\R\PBbb^{d-1}$
to some line $\ell$. Assume further that the sequence of tangent spaces $T_{y_k}S$, $k\in \NBbb$, converges to some 
subspace $\tau \subset T_xM$. Then $\ell \subset  \tau$.  
\end{enumerate}
By \cite[Lem.1.4.4]{PflAGSSS}, Whitney's condition (B) does not depend on the choice of the chart $\varphi$ around $x$.
A stratified subspace $(X,\calS)$ of a smooth manifold $M$ is said to be (A) respectively (B) regular if every pair of strata $(R,S)$ with 
$R$ incident to $S$ is (A) respectively (B) regular at each point $x \in R$.   
(B) regularity implies (A) regularity but in general not vice versa. Complex algebraic varieties \cite{WhiTAV}, orbit spaces of compact Lie group
actions \cite{PflAGSSS} and of proper Lie groupoids \cite{PflPosTanGOSPLG}, analytic varieties \cite{LojESA}, and subanalytic sets \cite{BieMilSSS} all
possess  (B) regular stratifications.